 \documentclass[11pt]{amsart}

\usepackage{palatino}
\usepackage[all]{xy}
\usepackage{amssymb}
\usepackage{latexsym}
\usepackage{amscd}
\usepackage{color}
\input amssym.def
\input amssym

\makeatletter \@addtoreset{equation}{section}

\makeatother \makeatletter

\def \<{\langle}
\def \>{\rangle}

\def \a{\alpha }

\def \l{\lambda }

\def \T{\mathcal T}
\def \b{\beta }

\newtheorem{theorem}{Theorem}[section]
\newtheorem{corollary}{Corrolary}[section]
\newtheorem{lemma}{Lemma}[section]
\newtheorem{conjecture}{Conjecture}[section]

\newtheorem{remark}{Remark}[section]

\newtheorem{definition}{Definition}[section]
\newtheorem{proposition}{Proposition}[section]

\newcommand{\bea}{\begin{eqnarray}}
\newcommand{\eea}{\end{eqnarray}}
\newcommand{\be}{\begin {equation}}
\newcommand{\ee}{\end{equation}}

\def \ga{\gamma }

\newcommand{\h}{\frak h}
\newcommand{\wt}{{\rm {wt} }   }

\newcommand{\Z}{\Bbb Z}
\newcommand{\W}{\mathcal W}

\newcommand{\Zp}{{\Bbb Z}_{>0} }

\newcommand{\N}{{\Bbb Z}_{\ge 0} }
\newcommand{\C}{\Bbb C}

\newcommand{\vak}{\bf 1}

\newcommand{\la}{\langle}
\newcommand{\ra}{\rangle}

\newcommand{\NS}{\frak{ns} }
\newcommand{\triplet}{\mathcal{W}(p)}
\newcommand{\striplet}{\mathcal{SW}(m)}
\newcommand{\hf}{\mbox{$\frac{1}{2}$}}
\newcommand{\thf}{\mbox{$\frac{3}{2}$}}

\begin{document}

\title[]{The $N=1$ triplet vertex operator superalgebras}

\author{ Dra\v zen Adamovi\' c  and Antun Milas}

\date{}

\address{Department of Mathematics, University of Zagreb, Croatia}
\email{adamovic@math.hr}

\address{Department of Mathematics and Statistics,
University at Albany (SUNY), Albany, NY 12222}
\email{amilas@math.albany.edu}
\thanks{The second author was partially supported by NSF grant DMS-0802962.}

 \markboth{Dra\v zen Adamovi\' c and Antun Milas} {Vertex
superalgebras}
\bibliographystyle{amsalpha}
  \maketitle

\begin{abstract}
We introduce a new family of $C_2$-cofinite $N=1$ vertex operator
superalgebras $\striplet$, $m \geq 1$, which are natural super
analogs of the triplet vertex algebra family $\triplet$, $p \geq 2$,
important in logarithmic conformal field theory. We classify
irreducible $\striplet$-modules and discuss logarithmic modules. We
also compute bosonic and fermionic formulas of irreducible
$\striplet$ characters. Finally, we contemplate possible connections
between the category of $\striplet$-modules and the category of
modules for the quantum group $U^{small}_q(sl_2)$, $q=e^{\frac{2 \pi
i}{2m+1}}$, by focusing primarily on properties of characters and the
Zhu's algebra $A(\striplet)$. This paper is a continuation of our
paper  Adv. Math. 217 (2008), no.6, 2664-2699.

\end{abstract}
\tableofcontents

\section{ Introduction}

Compared to rational vertex algebras, significantly less is known
about the structure of modules for general vertex algebras.
Recently, geared up with clues from the physics literature,  some
breakthrough has been achieved in understanding at least
quasi-rational vertex algebras (i.e., $C_2$-cofinite irrational
vertex algebras), and in particular the triplet vertex algebras
$\mathcal{W}(p)$, $ p \geq 2$ \cite{AdM-triplet}, \cite{FHST},
\cite{CF} (cf. also \cite{Abe} for $p=2$). But apart from the
symplectic fermions $\mathcal{W}(2)$, the description of
categories of weak (logarithmic) modules for other triplets
$\mathcal{W}(p)$, $p \geq 3$ remains open, even though there is
strong evidence for Kazhdan-Lusztig correspondence between the
category of logarithmic $\mathcal{W}(p)$-modules and certain
categories of modules for quantum groups (for these and related
developments we refer the reader to \cite{FHST}, and especially
\cite{FGST}, \cite{FGST1},  \cite{Se}, and references therein).

In \cite{AdM-triplet} we obtained several useful results about the
structure of the category of $\mathcal{W}(p)$-modules by using
primarily Zhu's algebra and Miyamoto's pseudocharacters \cite{Miy}.
Eventually, we will require more-or-less explicit
knowledge of "higher" Zhu's algebras for the triplet. But several
obstacles (e.g., explicit realization of certain logarithmic
modules) prevents us for taking this theory to the next level.
We hope that this approach, in particular,  will give an additional
evidence for Kazhdan-Lusztig correspondence, because
believe that proper understanding of the relationship between
quantum groups and triplets.

As with other familiar rational vertex operator algebras (e.g. Virasoro
minimal models), one may also wonder if the triplet has
(interesting!) $N=k$ super extensions, and whether those exhibit
similar properties (e.g., $C_2$-cofiniteness). In this paper we
solve this problem for $k=1$, by constructing a family of $N=1$
vertex operator superalgebras $\mathcal{SW}(m)$, $m \geq 1$, which
share many similarities with the triplet family.

In what follows we briefly recall the construction of $\mathcal{SW}(m)$ and present our main results.

Let us recall that the triplet vertex algebra
$\triplet$ \cite{FHST}, \cite{AdM-triplet} is defined as the kernel of a screening operator
acting from $V_L$ to $V_{L-\frac{\alpha}{p}}$
where $V_L$ is the vertex algebra associated to rank one even
lattice $\mathbb{Z}\alpha$, $\langle \alpha,\alpha \rangle=2p$, and $V_{L-\frac{\alpha}{p}}$ is a certain $V_L$-module. To
construct an $N=1$ super triplet we replace the even lattice with
an odd lattice such that  $\langle \alpha,\alpha \rangle=2m+1$, so that
$V_L$ has a natural vertex operator superalgebra structure. Then
we tensor $V_L$ with $F$, the free fermion vertex operator
superalgebra (cf. \cite{KWn}), and again, there is a screening operator
$$\tilde{Q} : \ V_L \otimes F  \longrightarrow V_{L-\frac{\alpha}{2m+1}} \otimes F.$$  The kernel of this
operator, denoted by $\striplet$,  is what we call  $N=1$ triplet vertex operator
superalgebra (or simply, $N=1$ super triplet). If we
restrict the kernel of $\tilde{Q}$ on the charge zero subspace we obtain another
vertex operator superalgebra
$$\overline{SM(1)} \subset \striplet,$$ called $N=1$ singlet vertex operator superalgebra. Both
vertex operator superalgebras contain Neveu-Schwarz vector $\tau$,
giving a representation of $\NS$ Lie superalgebra of central
charge $\frac{3}{2}-\frac{12 m^2}{2m+1}$. This is precisely the
central charge of $(1,2m+1)$ Neveu-Schwarz (degenerate) minimal
modules. A different $N=1$ extension of the symplectic fermion
$\mathcal{W}$-algebra $\mathcal{W}(2)$ was considered in \cite{MS}

By using the notation used by physicists, our super triplet would be an
example of a super $\mathcal{W}$-algebra of type
$\mathcal{W}(\frac{3}{2},2m+\frac{1}{2},2m+\frac{1}{2},2m+\frac{1}{2})$.
Similarly, the $N=1$ singlet algebra $\overline{SM(1)}$ is an
example of a super $\mathcal{W}$-algebra of type $\mathcal{W}(\frac{3}{2},2m+\frac{1}{2})$. We should say that for low $m$ (e.g., $m=1$)
some general properties of $\mathcal{W}$ superalgebras with two
generators were also discussed in the physics literature, but mostly by
using Jacobi identity and methods of Lie algebras (cf. \cite{BS} and references therein).
We should also mention that several general results about
$\mathcal{W}$-superalgebras associated to affine superalgebras were recently
obtained in \cite{Ar}, \cite{KWak} (see also \cite{HK}). However, super singlet and super triplet
vertex superalgebras do not appear in these works.

Because of the similarity between $\triplet$ and $\striplet$ many
results we obtain here are intimately related to those for the
triplet \cite{AdM-triplet} (cf. \cite{FGST}, \cite{CF}), but there are some subtle differences
which we address at various stages. However, to keep the paper
self-contained at many places we gave proofs that are almost
identical to those in \cite{AdM-triplet}.

Let us first consider the super singlet  $\overline{SM(1)}$.  This vertex superalgebra is too small to be $C_2$-cofinite (let alone rational!), which is evident from the following result.
\begin{theorem} Assume that $m \ge 1$.
\begin{itemize} \item[(i)] The singlet vertex superalgebra $\overline{SM(1)}$ is a simple $N=1$
vertex operator superalgebra generated by $\tau$ and a primary
vector $H$ of conformal weight $2m +\tfrac{1}{2}$.
 \item[(ii)]
The associative Zhu's algebra $A(\overline{SM(1)})$ is isomorphic to the
commutative algebra ${\C}[x,y] / \la P(x,y) \ra$
where $\la P(x,y) \ra $ is the ideal in ${\C}[x,y]$ generated by
the polynomial
$$ P(x,y) = y^2 - C_m \prod_{i=0} ^{2m} (x- h^{2i+1,1}), ,$$
where $C_m$ is a non-trivial constant and
%
$ h^{2i+1,1} = \frac{i (i- 2m)}{2 (2m+1)}. $
\end{itemize}
\end{theorem}

So the structure and representation theory of $\overline{SM(1)}$
is quite similar to that of $\overline{M(1)}$ investigated in
\cite{A-2003} and \cite{AdM-2007}. In particular we can construct
interesting logarithmic $\overline{SM(1)}$--modules and
{\em logarithmic intertwining  operators} as defined in \cite{HLZ}.

Next we study the vertex operator superalgebra $\striplet$. The main result on the structure on this vertex superalgebra is

\begin{theorem} \label{intro-1}   Assume that $m \ge 1$.
\begin{itemize} \item[(i)]  $\striplet$
is a simple $N=1$ vertex operator superalgebra generated by $\tau$
and three primary vectors $E,F, H$ of conformal weight $2m +
\tfrac{1}{2}$. \item[(ii)] The vertex operator superalgebra
$\striplet$ is irrational and $C_2$-cofinite.
   \item[(iii)]  $\striplet$ has precisely $2m+1$
inequivalent irreducible modules.
\end{itemize}
\end{theorem}
Our proof of (iii) imitates the proof of $C_2$-cofiniteness for the triplet $\mathcal{W}(p)$   \cite{AdM-triplet} (for a different proof see \cite{CF}).
The rest is done by combining methods of Zhu's associative
algebra and our knowledge of irreducible $V_L \otimes F$-modules. In parallel
with the triplet vertex algebra, we do not have an explicit
description of $A(\striplet)$, but we believe that the
following conjecture should hold true.
\begin{conjecture} \label{intro-2}
The Zhu's algebra decomposes as a sum of ideals
$$A(\striplet)=\bigoplus_{i=2m+1}^{3m}  \mathbb{M}_{h^{2i+1,1}} \oplus \bigoplus_{i=0}^{m-1} \mathbb{I}_{h^{2i+1,1}}
\oplus \mathbb{C}_{h^{2m+1,1}},$$ where $\mathbb{M}_{h^{2i+1,1}} \cong
M_2(\mathbb{C})$, ${\rm dim}(\mathbb{I}_{h^{2i+1,1}})= 2$ and
$\mathbb{C}_{h^{2m+1,1}}$ is one-dimensional. In particular
$A(\striplet)$ is $6m+1$-dimensional.
\end{conjecture}
In view of the classification result (cf. Theorem \ref{intro-2}),
it is important to compute irreducible characters and study
their modular transformation properties. As with the triplet
vertex algebra \cite{Fl95},  irreducible $\striplet$-characters
are often expressible as sums of modular forms of unequal weight.
Also, because we are working within vertex operator
superalgebras the $SL(2,\Z)$ group should be replaced with
the $\theta$-group $\Gamma_{\theta}$. Then we have
\begin{theorem} \label{intro-3} The $\Gamma_\theta$-closure of the space spanned by irreducible
$\striplet$-characters is $3m+1$-dimensional.
\end{theorem}
For a more precise statement see Theorem \ref{modular-inv}. Our result should be compared with \cite{Fl95}, where it was observed that the $SL(2,\mathbb{Z})$ closure
of the vector space of $\mathcal{W}(p)$ characters is $3p-1$ dimensional.
Finally, in parallel with \cite{FGK} and \cite{FFT}, we also obtain (see
Section 13) fermionic formulas for characters of irreducible
$\striplet$-modules. Our main results indicate that there is an
interesting relationship between characters of irreducible
$\mathcal{SW}(m)$-modules and irreducible characters of
$\mathcal{W}(2m+1)$-modules. It is not clear if there is a deeper
connection between these two $\mathcal{W}$-algebras.

Notice that if Conjecture \ref{intro-2} is true, then the center of
$A(\striplet)$ is  $3m+1$-dimensional, which is precisely the
dimension of the center of the small quantum group ${U^{small}_q
(sl_2)}$, $q = e^{\frac{2 \pi i }{2m+1}}$ \cite{Ker}. It is no
accident that this dimension matches the dimension in Theorem
\ref{intro-3} (similar phenomena occurs for the triplet vertex
algebra \cite{FGST}). Furthermore,  both $U^{small}_q(sl_2)$ and
$\striplet$ have the same number of irreducible modules \cite{Ker}
(see also \cite{La}). Thus, motivated by conjectures in \cite{FGST},
we expect the following (rather bold) conjecture to be true.
\begin{conjecture} \label{conj-strip} The category of weak $\striplet$-modules is
equivalent to the category of modules for the quantum group
${{U}^{small}_q(sl_2)}$, where $q=e^{\frac{2 \pi i}{2m+1}}$.
\end{conjecture}

\noindent {\bf   Acknowledgement:} We thank the anonymous referee for his/her valuable comments.

\section{Preliminaries}

In this section we briefly discuss the definition of vertex operator
superalgebras, their modules and intertwining operators as developed
in \cite{FFR}, \cite{K}, \cite{KWn}, \cite{Li}, \cite{HLZ}, \cite{HM}, etc. We
assume the reader is familiar with basics of vertex algebra
theory (cf. \cite{FHL}, \cite{FLM}, \cite{FB},  \cite{K}, \cite{LL}, etc.).

Let $V= V _{\bar 0} \oplus V _{\bar 1}$ be any ${\Z}_2$--graded
vector space. Then any element $u \in V _{\bar 0} $ (resp. $u \in V
_{\bar 1})$ is said to be even (resp. odd). We define $\vert u \vert
= {\bar 0} $ if $u$ is even and $\vert u \vert = {\bar 1}$ if $u$ is
odd. Elements in  $V _{\bar 0}$ or  $V _{\bar 1}$ are called
homogeneous. Whenever $\vert u \vert $ is written, it is understood
that $u$ is homogeneous.

The notion of vertex operator superalgebra is a natural (and
straightforward) generalization of the notion of vertex algebra
where the vector space $V$ in the definition is assumed to be
$\mathbb{Z}_2$-graded, where the vertex operator map
$$Y( \cdot,z) : V \longrightarrow {\rm Hom}(V,V((z)), \ \
Y(u,z)=\sum_{n \in \mathbb{Z}} u_n z^{-n-1}$$ is compatible with the
$\mathbb{Z}_2$-grading, and where Jacobi identity for a pair of
homogeneous elements is adjusted with an appropriate sign.

A vertex superalgebra $V$ is called a vertex operator superalgebra
if there is  a special  element $\omega \in V _{\bar 0} $ (called
{\em conformal vector}) whose vertex operator we write in the form
$$Y(\omega,z)=\sum_{n\in {\Z}}\omega_n z^{-n-1}=\sum_{n\in {\Z}
}L(n)z^{-n-2},$$ such that $L(n)$ close the Virasoro algebra
representation on $V$, and where $V$ is
$\frac{1}{2}{\mathbb{Z}}$-graded (by {\em weight}), truncated from
below, with finite-dimensional vector spaces. Also, the grading is
determined with the action of the Virasoro operator $L(0)$.
In this paper, we shall assume that
$$ V_{\bar 0}=\coprod_{ n \in {\N} } V(n), \ \   V_{\bar 1}=\coprod_{ n \in {\tfrac{1}{2}}+{\N} } V(n), \ \ \mbox{where} \ \ V(n) = \{ a
\in V \ \vert \ L(0) a = n a \}.
$$
For $a \in V(n)$, we shall write $\wt (a) = n$ or ${\rm
deg}(a)=n$.
We shall sometimes refer to vertex operator superalgebra $V$ as a
quadruple $(V, Y, {\bf 1}, \omega)$, where ${\bf 1}$ is the vacuum
vector (as for vertex operator algebras).

We say that the vertex operator superalgebra $V$ is generated by
the set $S \subset V$ if
$$ V=\mbox{span}_{\C} \{ u^{1} _{n_1} \cdots u^{r} _{n_r}
{\vak} \ \vert \  \ u^{1}, \dots, u^{r} \in S, \ n_1, \dots, n_r
\in {\Z}, r \in {\Zp} \}.$$

The vertex operator algebra  $V$  is said to be {\em strongly
generated} (cf. \cite{K}) by the set $R$ if
$$V = \mbox{span}_{\C} \{ u^{1} _{n_1} \cdots u^{r} _{n_r}
{\vak} \ \vert \  \ u^{1}, \dots, u^{r} \in R, \ n_i < 0, r \in
{\Zp} \}.$$

In parallel with vertex algebras we can define the notion of {\em weak
module} for vertex operator superalgebras. Again, the only new
requirement is that the vector space $M$ in the definition is
$\mathbb{Z}_2$-graded, with grading compatible with respect to the
action of $V$, and where the Jacobi identity is adjusted as in the
case of vertex superalgebras. The vertex operator acting on $M$ is
usually denoted by $Y_M$.
%
%
%

A weak $V$--module $(M,Y_M)$ is called an (ordinary) $V$-module if
$M$ carries an action of the Virasoro algebra via the expansion of
$Y_M(\omega,x)$, and in addition $M$ is equipped with a
$\mathbb{R}$-grading (or even $\mathbb{C}$-grading) determined by
the Virasoro operator $L(0)$. In addition, the grading is truncated
from below, with finite dimensional graded subspaces.

As usual, we say that a $V$-module $M$ is irreducible (or simple) if $M$ has
no proper submodules. We say that a vertex operator superalgebra is
{\em rational} if every $V$-module $M$ is semisimple (i.e., $M$
decomposes as a direct sum of irreducible modules) and if $V$ has
only finitely many (inequivalent) irreducible modules.

\begin{definition} {\em Let $V$ be a vertex operator superalgebra. We say that a weak $V$-module $M$ is
{\em logarithmic}, if it carries an action of the Virasoro algebra
and if it admits decomposition
$$M=\coprod_{r \in \mathbb{C}} M_r,$$
where
$$M_r=\{ v \ :  \ (L(0)-r)^k v=0, \ {\rm for \ some} \ k \in
\mathbb{N} \}.$$}
\end{definition}

\subsection{Zhu's algebra $A(V)$}

We define two bilinear maps $* : V  \times V \rightarrow V$,
$\circ : V \times V \rightarrow V$ as follows. For homogeneous $a,
b \in V$ let
\bea
a* b &&= \left\{\begin{array}{cc}
 \  \mbox{Res}_x Y(a,x) \frac{(1+x) ^{\deg (a)}}{x}b  & \mbox{if} \ a,b  \in V_{\bar{0}} \\
  0 & \mbox{if} \ a \ \mbox{or} \ b  \in V_{\bar{1}} \
\end{array}
\right.  \\
a\circ b &&= \left\{\begin{array}{cc}
 \  \mbox{Res}_x Y(a,x) \frac{(1+x) ^{\deg (a)} }{x^2}b  & \mbox{if} \ a  \in V_{\bar{0}} \\
  \  \mbox{Res}_x Y(a,x) \frac{(1+x) ^{\deg (a) -\tfrac{1}{2} }}{x}b  & \mbox{if} \ a   \in V_{\bar{1}} \
\end{array}
\right. \eea

Next, we extend $*$ and $\circ$ on $V \otimes V$ linearly, and
denote by $O(V)\subset V$ the linear span of elements of the form $a
\circ b$, and by $A(V)$ the quotient space $V / O(V)$. The image of
$v \in V$, under the natural map $V \mapsto A(V)$ will be denoted by
$[v]$.  The space $A(V)$ has a unital  associative algebra
structure, with the product  $*$ and  $[\bf 1]$ as the unit element.
The associative algebra $A(V)$ is called  the Zhu's algebra of $V$.

Assume that $M= \oplus_{n \in \tfrac{1}{2}{\N}} M(n)$ is a
$\tfrac{1}{2}{\N}$--graded $V$--module. Then the top component
$M(0)$ of $M$ is a $A(V)$--module under the action $[a] \mapsto o(a)
=a_{\wt (a)-1}$ for homogeneous $a$ in $V_{\bar{0}}$.
 We shall
sometimes write $a(0)$ for $o(a)$.
(Note that if   $ a \in  V_{\bar{1}} $, then $[a] = 0$ in $A(V)$. We
 formally set  $o(a) = a(0) = 0$ in this case.)

Moreover, there is one-to-one correspondence between irreducible
$A(V)$--modules and irreducible $\tfrac{1}{2}\N$--graded
$V$--modules (cf. \cite{KWn}).

As usual, for a vertex operator superalgebra $V$ we let
$$C_2(V)=\{ a_{-2}b \ : \ a, b \in V \}.$$
Then it is not hard to see that $$\mathcal{P}(V)=V/C_2(V)$$ has a
super Poisson algebra structure with the multiplication
$$\bar{a} \cdot \bar{b}=\overline{a_{-1} b},$$
and the Lie bracket
$$[\bar{a},\bar{b}]=\overline{a_0 b},$$
where $-$ denotes the natural projection from $V$ to
$\mathcal{P}(V)$ (see for instance \cite{Z}). Therefore we have a
decomposition $\mathcal{P}(V)=\mathcal{P}(V)_0 \oplus
\mathcal{P}(V)_1$ into even and odd subspace, respectively. If
$V/C_2(V)$ is finite-dimensional we say that $V$ is $C_2$-cofinite.
%
%
Let $a, b \in V$, be ${\Z}_2$ homogeneous. Then by using
super-commutator formulae in vertex operator superalgebras one can
easily see that
\bea && \bar{a} \cdot \bar{b} - (-1) ^ { \vert a \vert \vert b \vert
} \bar{b} \cdot \bar{a}  = 0 \quad \mbox{in} \ V / C_2(V).
\label{super-commute} \eea

 The  following result was proved in \cite{dSK}, and it is a
generalization of Proposition 2.2 in  \cite{Abe}.

\begin{proposition} \label{abe}
Let $V$ be strongly generated by the set $S$. Then we have:

\item[(1)]  $\mathcal{P}(V)$ is generated by the set $\{ \overline{a},
a \in S \}$.

\item[(2)]  $A(V)$ is generated by the set $\{ [a], a \in S \}$.

\item[(3)] If  $V$ is $C_2$-cofinite
$${\rm dim}(\mathcal{P}(V)_0) \geq {\rm dim}(A(V)).$$
\end{proposition}

\subsection{Intertwining operators among vertex operator
superalgebra modules}

Intertwining operators for superconformal vertex operator algebras
were introduced in \cite{KWn}. Their theory is further developed
in \cite{HM} by using both even and odd formal variables. We
briefly outline the definition here.

\begin{definition}
Let $V$ be a vertex operator superalgebra and $M_1$, $M_2$ and
$M_3$ a triple of $V$--module. An intertwining operator
$\mathcal{Y}( \cdot , z)$ of type ${M_3 \choose M_1 \ M_2 }$ is a
linear map
$$ \mathcal{Y} : M_1 \rightarrow \mbox{End} (M_2,M_3) \{ z \},$$
$$ w_1 \mapsto \mathcal{Y}(w_1,z) = \sum_{ n \in {\C}} (w_1)_n z
^{-n-1}, $$
satisfying the following conditions for $w_i \in M_i$, $i=1,2$ and
$a \in V$:

\begin{enumerate}

\item[(I1)]$\mathcal{Y}(L(-1)w_1,z)=\frac{d}{dz}\mathcal{Y}(w_1,z)$.
\item[(I2)] $(w_1)_n(w_2) = 0$ for ${\rm Re}(n)$ sufficiently
large.
\item[(I3)] The following Jacobi identity holds
\begin{eqnarray}
&
&z_{0}^{-1}\delta\left(\frac{z_{1}-z_{2}}{z_{0}}\right)Y_{M_3}(a,z_{1})\mathcal{Y}(w_1,z_{2})w_2
 -(-1) ^{\vert a \vert \vert w_1 \vert}
z_{0}^{-1}\delta\left(\frac{z_{2}-z_{1}}{-z_{0}}\right)
\mathcal{Y}(w_1,z_{2})Y_{M_2}(a,z_{1})w_2 \nonumber \\
&=&z_{2}^{-1}\delta\left(\frac{z_{1}-z_{0}}{z_{2}}\right)\mathcal{Y}(Y_{M_1}(a,z_{0})w_1,z_{2})w_2,
\nonumber
\end{eqnarray}
for $\mathbb{Z}_2$-homogeneous $a$ and $w_1$.

\end{enumerate}
\end{definition}

We shall denote by
$$I \ {M_3 \choose M_1 \ M_2 }$$
the vector space of intertwining operators of type ${M_3 \choose
M_1 \ M_2 }$. Their dimensions are known as  the "fusion rules".

\section{$N=1$ Neveu-Schwarz vertex operator superalgebras}

The $N=1$ Neveu-Schwarz (or simply NS) algebra is the Lie
superalgebra
$${\NS} = \bigoplus_{n \in {\Z} } {\C} L(n) \bigoplus \bigoplus_{m
\in {\hf} + {\Z} }  {\C}G(m) \bigoplus {\C}C$$
 with commutation relations ($m,n \in {\Z} $):
\bea &&[L(m),L(n)] =(m-n)L({m+n})+\delta_{m+n,0}\frac{m^3-m}{12}C,
\nonumber  \\ && \label{ns1com} [G(m+\frac{1}{2}),L(n)]   =
(m+\frac{1}{2}-\frac{n}{2})G({m+n+\frac{1}{2}}),    \\ &&
\label{ns2com} \{ G({m+\frac{1}{2}}),G({n-\frac{1}{2}}) \}  = 2
L({m+n})+\frac{1}{3}m(m+1)\delta_{m+n,0}C,  \\ && [L(m),C] = 0,
\,\,\,\,\, [G({m+\frac{1}{2}}),C] = 0. \nonumber \eea

It is important to consider vertex algebras which admit an action
of the $N=1$ Neveu-Schwarz algebras (cf. \cite{HM}). These vertex
operator superalgebras are called $N=1$ Neveu-Schwarz vertex
operator superalgebras and are subject to an additional axiom:

There exists $\tau \in V_{3/2}$ (superconformal
vector) such that
$$Y(\tau,z)=\sum_{n \in \mathbb{Z}+1/2} G(n)z^{-n-3/2}, \ \ G(n) \in {\rm End}(V)$$
where $G(n)$ satisfy bracket relations as in (\ref{ns1com}) and
(\ref{ns2com}).

The simplest examples of $N=1$ vertex operator superalgebras are
$\NS$-modules $L^{\NS}(c,0)$, $c \neq 0$ where we use the standard
notation and for any $(c,h) \in {\C} ^{2}$ we denote by
$L^{\NS}(c,h)$ the corresponding irreducible highest weight
${\NS}$--module with central charge $c$ and highest weight $h$
(cf. \cite{KWn}, \cite{Li}, \cite{A-1997}, \cite{HM}).  It is
well-known that the vertex operator superalgebra $L^{\NS} (c,0)$,
$c \neq 0$ is simple.

Set $$ c_{p,q} = \frac{3}{2}(1-\frac{2(p-q)^2}{pq}), $$ $$
h_{p,q}^{r,s} = \frac{(sp-rq)^2 -(p-q)^2}{8pq}. $$ In the rest of
the paper we shall focus on certain $\NS$ modules of central
charge $c_{2m+1,1}$, $m \geq 1$.

\section{Fusion rules for $N=1$ superconformal $(2m+1,1)$-models}

From now on we will mostly focus on (non-minimal) $(2m+1,1)$-models,
so that $p=2m+1$, $q=1$. Relevant lowest weights are
$h^{r,s}:=h_{2m+1,1}^{r,s}$, $r,s \in \mathbb{Z}$.

It will be of great use to determine the fusion rules \be
\label{fusion-full} I \ { L(c_{2m+1,1},h^{r'',s''}) \choose
L(c_{2m+1,1},h^{r,s}) \ L(c_{2m+1,1},h^{r',s'})} \ee for certain
triples $(r,s)$, $(r',s')$ and $(r'',s'') \in \mathbb{Z}^2$. For
$m=0$ (i.e., the $c=3/2$ case) these numbers were computed in (see
\cite{M0}). In particular, for every $s>0$ we have: \be \label{c32}
L(\frac{3}{2},h^{1,3}) \times
L(\frac{3}{2},h^{1,2s+1})=L(\frac{3}{2},h^{1,2s-1}) \oplus
L(\frac{3}{2},h^{1,2s+1}) \oplus L(\frac{3}{2},h^{1,2s+3}), \ee
where $\times$ is just  a formal product indicating which triples of irreducible modules admit nontrivial
fusion rules (all with multiplicity one). As shown in \cite{M0}, the
fusion rules for $m=0$ can be computed by using certain projection
formulas for singular vectors combined with Frenkel-Zhu's formula.
It is not hard to see that the same approach extends to $m \geq 1$
as well.  We only have to apply appropriate projection formulas as
in Lemma 3.1 of \cite{IK2}. Actually, for purposes of this paper we
do not need any of results from \cite{IK2}, because we are
interested only in special properties of "fusion rules"
(\ref{fusion-full}) (nevertheless, see Remark \ref{iohara-remark}).
\begin{proposition} \label{fusion}
For every $i=0,...,m-1$ and $n \geq 1$ we have: the space
$$I \ {L(c_{2m+1,1},h) \choose L(c_{2m+1,1},h^{1,3}) \
L(c_{2m+1,1},h^{2i+1,2n+1})}$$ is nontrivial only if $h \in
\{h^{2i+1,2n-1}, h^{2i+1,2n+1},h^{2i+1,2n+3} \}$, and  $$I \
{L(c_{2m+1,1},h)  \choose  L(c_{2m+1,1},h^{1,3}) \
L(c_{2m+1,1},h^{2i+1,1})}$$ is nontrivial only if $h=h^{2i+1,3}.$


Similarly, for every $i=0,...,m-1$ and $n \geq 2$ we have: the space
$$I \ {L(c_{2m+1,1},h) \choose L(c_{2m+1,1},h^{1,3}) \
L(c_{2m+1,1},h^{2i+1,-2n+1})}$$ is nontrivial only if $h \in \{
h^{2i+1,-2n-1},h^{2i+1,-2n+1}, h^{2i+1,-2n+3} \}$, and
$$I \ {L(c_{2m+1,1},h) \choose L(c_{2m+1,1},h^{1,3}) \
L(c_{2m+1,1},h^{2i+1,-1})}$$ is nontrivial only if $h \in
\{h^{2i+1,-3},h^{2i+1,-1} \}$.

%
%
For a stronger statement see Remark \ref{iohara-remark}.

\end{proposition}

\noindent {\em Proof.} We assume that $n \geq 1$ (for other cases
essentially the same argument works). Let $A(L(c_{2m+1,0},0))$ be
the Zhu's algebra of $L(c_{2m+1,0},0)$ (polynomial algebra in one
variable) and $A(L(c_{2m+1},h))$ the $A(L(c_{2m+1,0},0))$-bimodule
of $L(c_{2m+1},h)$ \cite{FZ}.

As in \cite{M0}, it is sufficient to analyze the structure of the
$A(L(c_{2m+1,0},0))$-module \be \label{bim-tensor}
A(L(c_{2m+1},h^{1,3})) \otimes_{A(L(c_{2m+1,0},0)}
L(c_{2m+1,1},h^{2i+1,2n+1})(0), \ee where $L(c_{2m+1,1},h)(0)$
denotes the top weight component of $L(c_{2m+1,1},h)$ (cf.
\cite{FZ}). From \cite{IK2} (or elsewhere) it follows that the Verma
module $M(c_{2m+1,0},h^{1,3})$ combines in the following short exact
sequence
$$0 \longrightarrow M(c_{2m+1,0},h^{1,3}+\frac{3}{2}) \longrightarrow  M(c_{2m+1,0},h^{1,3}) \longrightarrow L(c_{2m+1,0},h^{1,3}) \longrightarrow 0.$$
Thus the maximal submodule of $M(c_{2m+1,0},h^{1,3})$ is generated by a singular vector of
weight $h^{1,3}+\frac{3}{2}$ (explicitly,
$(-L(-1)G(-1/2)+(2m+1)G(-3/2))v_{1,3}$
where $v_{1,3}$ is the highest weight vector in
$M(c_{2m+1,0},h^{1,3})$). Now, as in \cite{M0}, it is not hard to
see that the space (\ref{bim-tensor}) is three-dimensional and that
all fusion rules covered by the statements are at most $1$
(actually, they are all one; see Remark \ref{iohara-remark}).

\qed

\begin{remark} \label{iohara-remark}
{\em We can actually prove "if and only if" statement in Proposition
\ref{fusion} by using at least two different methods. On one hand we
would have to combine methods from \cite{M0} and projection formula
in Lemma 3.1 \cite{IK2} (we do not have explicit singular vectors to
work with!). Alternatively, with Proposition \ref{fusion}, it is
sufficient to construct non-trivial intertwining operators for all
types covered in Proposition \ref{fusion}. This was actually done in
later sections.

We should say that our fusion rules formulas coincide with
Iohara-Koga's fusion rule formula in the generic case, which are
computed by using coinvariants and projection formulas rather than
Frenkel-Zhu's formula \cite{IK2}. But as we know the coinvariant
approach and Frenkel-Zhu's formulas yield the same answer in
practically all known examples (for further examples see \cite{W1},
\cite{M0}, \cite{M4}). }
\end{remark}

\section{Lattice and fermionic vertex superalgebras }
\label{stwist}


 We shall first recall some basic facts about lattice and
 fermionic vertex superalgebras.

 Let $m \in {\N}$. Let $\widetilde{L}={\Z}\beta$  be a
rational lattice of rank one with nondegenerate
 bilinear
 form $\la \cdot, \cdot \ra$ given by
 $$ \la \beta , \beta \ra = \frac{1}{ 2 m + 1}.$$
 Let
 ${\h}={\C}\otimes_{\Z} \widetilde{L}$.
  Extend the form $ \la \cdot,
\cdot \ra $ on $\widetilde{L}$ to ${\h}$.
 Let $\hat{{\h}}={\C}[t,t^{-1}]\otimes {\h} \oplus {\C}c$ be the affinization of
${\h}.$
Set
$
\hat{{\h}}^{+}=t{\C}[t]\otimes
{\h};\;\;\hat{{\h}}^{-}=t^{-1}{\C}[t^{-1}]\otimes {\h}.
$
Then $\hat{{\h}}^{+}$ and $\hat{{\h}}^{-}$ are abelian subalgebras
of $\hat{{\h}}$. Let $U(\hat{{\h}}^{-})=S(\hat{{\h}}^{-})$ be the
universal enveloping algebra of $\hat{{\h}}^{-}$. Let ${\l} \in
{\h}$. Consider the induced $\hat{{\h}}$-module
\begin{eqnarray*}
M(1,{\l})=U(\hat{{\h}})\otimes _{U({\C}[t]\otimes {\h}\oplus
{\C}c)}{\C}_{\l}\simeq
S(\hat{{\h}}^{-})\;\;\mbox{(linearly)},\end{eqnarray*} where
$t{\C}[t]\otimes {\h}$ acts trivially on ${\C}_{\l} \cong \C$,
${\h}$ acting as $\la h, {\l} \ra$ for $h \in {\h}$
and $c$ acts on ${\C}_{\l}$ as multiplication by 1. We shall write
$M(1)$ for $M(1,0)$.
 For $h\in {\h}$ and $n \in {\Z}$ write $h(n) =  t^{n} \otimes h$. Set
$
h(z)=\sum _{n\in {\Z}}h(n)z^{-n-1}.
$
Then $M(1)$ is a vertex   algebra which is generated by the fields
$h(z)$, $h \in {\h}$, and $M(1,{\l})$, for $\l \in {\h}$, are
irreducible modules for $M(1)$.

As in  \cite{DL} (see also \cite{D}, \cite{FLM}, \cite{GL},
\cite{K}), we have the generalized vertex algebra
$$ V_{ \widetilde{L} } = M(1) \otimes {\C}[\widetilde{L}], $$
where ${\C}[\widetilde{L}]$ is a group algebra of $\widetilde{L}$
with a generator $e ^{\beta}$.
For $v \in  V_{ \widetilde{L} }$, let $Y(v,z) = {\sum_{ s \in
\frac{1}{2m+1} {\Z} } v_s z^{-s -1}}$ be the corresponding vertex
operator (for precise formulae see \cite{DL}).

Define $\alpha = (2 m +1)  \beta$. Then $\la \a , \a \ra = 2 m+1$,
implying  $L= {\Z} \alpha \subset \widetilde{L}$ is an integer
lattice. Therefore the subalgebra $V_L \subset V_{ \widetilde{L}
}$ has the structure of a vertex  superalgebra.

 Define the Schur polynomials $S_{r}(x_{1},x_{2},\cdots)$
 in variables $x_{1},x_{2},\cdots$ by the following equation:
\begin{eqnarray}\label{eschurd}
\exp \left(\sum_{n= 1}^{\infty}\frac{x_{n}}{n}y^{n}\right)
=\sum_{r=0}^{\infty}S_{r}(x_1,x_2,\cdots)y^{r}.
\end{eqnarray}

For any monomial $x_{1}^{n_{1}}x_{2}^{n_{2}}\cdots x_{r}^{n_{r}}$
we have an element $$h(-1)^{n_{1}}h(-2)^{n_{2}}\cdots
h(-r)^{n_{r}}{\vak} $$ in $M(1)$   for $h\in{\h}.$
 Then for any polynomial $f(x_{1},x_{2}, \cdots)$,  $f(h(-1),
h(-2),\cdots){\vak}$ is a well-defined element in $M(1)$ . In
particular, $S_{r}(h(-1),h(-2),\cdots){\vak} \in M(1)$  for $r \in
{\N}$. Set $S_r (h)$ for $S_{r}(h(-1),h(-2),\cdots){\vak}$.

%
%
%
%
%
%
%
%
%

 The following relations in the generalized vertex
operator algebra $V_{ \widetilde{L} } $ are of great importance:
\begin{eqnarray}\label{eab1}
e ^{\gamma} _{i} e ^{\delta}=0\;\;\;\mbox{ for }i\ge
-\<\gamma,\delta\>.
\end{eqnarray}
Especially, if $\<\gamma,\delta\>\ge 0$, we have $e ^{\gamma} _{i}
e ^{\delta} =0$ for  $i\in {\N}$, and if
$\<\gamma,\delta\>=-n<0$,
 we get
 \bea\label{eab} && e ^{\gamma}
_{i-1} e ^{\delta} =S_{n-i}(\gamma) e ^{\gamma + \delta}
\;\;\;\mbox{ for }i\in \{ 0, \dots, n\}. \eea

\subsection{Fermionic vertex operator superalgebra F} In what
follows we consider the Clifford algebra $CL$, generated by
$\{ \phi  (n) , n \in  {\hf} + {\Z} \} \cup \{1\}$ and relations
$$ \{ \phi  ({n}) , \phi (m)\}=   \delta_{n,-m},
\quad n,m \in {\hf}+\mathbb{Z} . $$

Let $F$ be the $CL $--module generated by the vector ${\vak}$ such
that
$$ \phi ({n}) {\vak} = 0, \ n > 0.$$
 Then the field
$$Y( \phi ({-\hf}) {\bf 1} ,z) =\phi (z) = \sum_{ n \in {\hf} +
{\Z} } \phi ({ n }) z ^{-n- \hf}, $$
  generate
a unique vertex operator superalgebra structure on $F$. We choose
$$ {\omega}^{(s)} = \frac{1}{2}  \phi ({-\thf}) \phi
({- \hf}) {\vak}$$
for the Virasoro element giving  central charge $\frac{1}{2}$. Moreover,
$F$ is a rational vertex operator superalgebra, and $F$ is up to
equivalence the unique irreducible $F$--module (see \cite{FRW},
\cite{KWn}, \cite{Li}).

\subsection{ Vertex superalgebra $SM(1)$}

In this subsection we study the vertex superalgebra $SM(1):=M(1)
\otimes F$. We shall first define a family of $N=1$ superconformal
vectors in $SM(1)$. For every $m \in {\N}$, we define (see also
\cite{MR}, \cite{K},\cite{IK})
%
%
\bea &&\tau = \frac{1}{ \sqrt{2 m +1 }} \left( \alpha(-1) {\vak}
\otimes \phi (-\hf) {\vak} + 2 m {\vak} \otimes  \phi ({-\thf})
{\vak}.
\right) \nonumber \\
&&
\omega = \frac{1}{2 (2m+1)} ( \alpha(-1) ^{2} +  2 m \alpha(-2) )
{\vak} \otimes {\vak} +    {\vak} \otimes  {\omega}^{(s)}.
\nonumber
 \eea
Set
$$Y(\tau,z)=G(z) = \sum _{n \in {\Z} } G(n+\hf) z^{-n-2} , \
\  Y(\omega,z)=L(z) = \sum_{n \in {\Z} } L(n) z^{-n-2} .$$ Then
$\tau$ is an $N=1$ superconformal vector, and the vertex subalgebra
of $SM(1)$ strongly  generated by the fields $G(z)$  and $L(z)$ is
isomorphic to the Neveu-Schwarz vertex operator superalgebra
$L^{\NS}(c_{2 m+1,1}, 0)$, where $c_{2 m +1,1} =
\frac{3}{2}(1-\frac{8 m^{2}}{2 m+1})$. In other words, $SM(1)$
becomes a Fock module for the Neveu-Schwarz algebra with central
charge $c_{2 m+1,1}$.  Moreover, for every $\l \in {\h}$, the
$SM(1)$--modules $SM(1,\l):=M(1,\l) \otimes F$ is also a Fock module
with central charge $c_{2m+1,1}$
%
%
and conformal weight
\be \frac{1}{2 (2m+1)} ( \la \l, \a  \ra ^2 - 2m \la \l , \a \ra )
\label{conf-weight} . \ee
 Now we want to describe the
structure of these Fock modules viewed as $\NS$-modules. For this
purpose we need the concept of screening operators. As in
\cite{A-2003}, we shall construct these operators using generalized
vertex algebras.

The $N=1$ superconformal vector  $\tau \in M(1) \otimes F$ also
defines an $N=1$ superconformal structure on $V_{\widetilde{L} }
\otimes F$ and $V_L \otimes F$. In particular,  $V_L \otimes F$ is
an $N=1$  vertex operator superalgebra.  The operator $L(0)$
defines a ${\hf} {\N}$--gradation on $V_L \otimes F$. Recall that
 $\mbox{wt} (v) = n$ if $L(0) v = n v$.

Define
%
\bea && s ^{(1)} = e^{\a} \otimes \phi(-\hf){\vak} \in V_L \otimes
F, \nonumber \\ && s ^{(2)} = e^{-\beta} \otimes \phi(-\hf){\vak}
\in V_{ \widetilde{L} } \otimes F . \nonumber   \eea

By  using the Jacobi identity in the (generalized) vertex algebras
$ V_L \otimes F$ and $V_{\widetilde{L}} \otimes F$ we get the
following formulas
\bea
&& [G(n+\frac{1}{2}), s^{(1)} _ i ] = -\frac{i}{\sqrt{2m+1}}
e^{\a}_{i+n}, \quad [L(n), s^{(1)} _ i] = - i \  s^{(1)} _ {i+n}
\quad (i \in {\Z})  \label{komut-1}  \\
&& [G(n+\frac{1}{2}), s^{(2)} _ r ] = r \sqrt{2m+1}  e^{-\b}_{i+n},
\quad [L(n), s^{(2)} _ r] = - r \  s^{(2)} _ {r+n} \quad ( r  \in
\frac{1}{2m+1 }{\Z} )  \label{komut-2} .
\eea
Let
 \bea
&& Q = s^{(1)} _0 = \mbox{Res}_z \ Y( s^{(1)},z) , \nonumber
\\ && \widetilde{Q} = s^{(2)}_0 = \mbox{Res}_z Y( s^{(2)},
z) \  . \nonumber \eea

From relations (\ref{komut-1}) and (\ref{komut-2}) we see that the
operators $Q$ and $\widetilde{Q}$ commute with the action of the
Neveu-Schwarz algebra (see also \cite{IK}).

We are interested in the action of these
 operators on $SM(1)$. In fact, $Q$ and $ \widetilde{Q}$ are the screening
 operators, and therefore  $\mbox{Ker}_{SM(1)} Q$ and
$\mbox{Ker}_{SM(1)} \widetilde{Q}$ are vertex subalgebras of $SM(1)$
(for details see Section 14 in \cite{FB} and reference therein).

The following lemma gives the basic  properties of the operators
$Q$ and $ \widetilde{Q}$. The proof is similar to that of Lemma
2.1 in \cite{A-2003}.

\begin{lemma} \label{rel-tilde}
%
%
\item[(i)] If $m\ne 0$, $[ Q , \widetilde{Q}] =0$.
\item[(ii)] $ \widetilde{Q} e ^{ n \alpha}  \ne 0$, $n \in
{\Zp}$. \item[(iii)]$ \widetilde{Q} e ^{ -n \alpha}  =0$, $n \in
{\N}$.
\end{lemma}

We now define the following three (non-zero) elements in the
vertex operator superalgebras $V_L \otimes F$:
 $$ F =e ^{-\a}, \ \ H =   Q F, \ \ E =  Q ^{2}  F. $$
By using expression for conformal weights (\ref{conf-weight}) and
Lemma \ref{rel-tilde}, we conclude that  these vectors are singular
vectors for the action of the Neveu-Schwarz algebra, and
$$ \mbox{wt} (F) = \mbox{wt} (H) = \mbox{wt} (E) = h^{1,3}=2 m + \hf. $$
It is also important to notice that $H \in SM(1)$.

The proof of the following result is similar to that of Lemma 3.1
in \cite{A-2003}.

 \begin{lemma} \label{pomoc1} In the vertex operator superalgebra $V_L \otimes F$   the
 following relations hold:
\item[(i)] $ Q ^{3} F = 0$.
 \item[(ii)]  $E_{i} E = F_{i} F = 0$, for every $i \ge - 2  m -1$.
 \item[(iii)]  $ Q  ( H_{i} H )= 0$, for every $ i \ge - 2m -1$.

\end{lemma}

We define
\bea \label{def-hat} \widehat{F} = e^{-\a} \otimes
\phi(-\tfrac{1}{2}), \quad \widehat{H}= Q \widehat{F}, \quad
\widehat{E} = Q^{2} \widehat{F}. \eea

These vectors are even and  have conformal weight $2m+1$. We will
need the following result.

\begin{lemma} \label{nilp-hat}
We have
$$\widehat{F}_i \widehat{F}=0, \ \ \widehat{E}_i \widehat{E}=0, \ \ i \geq -2m.$$
Also,
$$Q(\widehat{H}_i \widehat{H})=0, \ \ i \geq -2m.$$
\end{lemma}
\noindent {\em Proof.}
 Since $Q$
acts as a derivation if $\widehat{F}_i \widehat{F}=0$, for $i \geq
-2m$ then $Q^4(\widehat{F}_i \widehat{F})=6 \widehat{E}_i
\widehat{E}=0$, for $i \geq -2m$. We only have to notice relations
$$\widehat{F}_k \widehat{F}={\rm Res}_{x} x^k Y(e^{-\alpha},x)e^{-\alpha} \otimes
Y(\phi(-1/2),x)\phi(-1/2){\bf 1},$$ $${\rm Res}_x x^i
Y(e^{-\alpha},x)e^{-\alpha}=0, \ \ i \geq -2m-1,$$ proven in Lemma
\ref{pomoc1}, and
$${\rm Res}_x x^j Y(\phi(-1/2),x)\phi(-1/2){\bf 1}=0, \ \ j \geq 1.$$

The last formula follows from $Q^3(\widehat{F}_i \widehat{F})=0$
for $i \geq -2m$. \qed

\section{The $N=1$ Neveu-Schwarz module structure of $V_L \otimes F$-modules}

For $i \in {\Z}$, we  set
\bea && \ga_i = \frac{i}{2m+1} \a \label{ga-def}. \eea

We shall first present results on the structure of $V_L \otimes
F$--modules as modules for the $N=1$ Neveu-Schwarz algebra. It is
a known fact that irreducible $V_L \otimes F$-modules are given by
$$V_{L+\gamma_i} \otimes F, \ \  i=0,...,2m.$$
Each $V_{L+\gamma_i}$ is a direct sum of super Feigin-Fuchs
modules via
$$V_{L+\gamma_i} \otimes F=\bigoplus_{n \in \mathbb{Z}} (M(1) \otimes
e^{\ga_i +n \alpha}) \otimes F.$$

 We shall now investigate
  the action of the operator $Q$.
Since operators $Q ^{j}$, $j \in {\Zp}$, commute with the action
of the Neveu-Schwarz algebra, they are  actually intertwiners
between super Feigin-Fuchs  modules inside  $V_{L+\ga_i} \otimes
F$.  Assume that $0 \le i \le m$. If $Q^{j} e^{\ga_i-n \alpha}$ is
nontrivial, it is a singular vector in the Fock module $SM(1,
\ga_i + (j -n ) \alpha)$ of weight
$$\mbox{wt}(Q^{j} e^{\ga_i -n \alpha}) = \mbox{wt} (e^{\ga_i -n \alpha}) =h ^{2i+1,2n+1},  $$
where $h^{2i+1,2n+1}:= h_{1,2m+1} ^{2i+1,2n+1}$.
Since $\mbox{wt} ( e^{ \ga_i + (j-n) \alpha} ) > \mbox{wt} (
e^{\ga_i -n \alpha} )$ if $j  > 2n   $, we conclude that
\bea \label{trivial-1}
&& Q^ {j} e^{\ga_i-n \alpha} = 0 \quad \mbox{for} \ j > 2 n. \eea
One can similarly see that for $ m+1 \le i \le 2m$ : \bea
\label{trivial-2}
&& Q^ {j} e^{\ga_i-n \alpha} = 0 \quad \mbox{for} \ j > 2 n +1.
\eea

The following lemma is useful for constructing singular vectors in
$V_{L+\ga_i} \otimes F$:
\begin{lemma} \label{non-triv}
\item[(1)]  $ Q^{2n} e^{\ga_i - n\a} \ne 0$ \  for  \  $0 \le i
\le m$.
\item[(2)]$ Q^{2n+1} e^{\ga_i - n\a} \ne 0$ \ for \  $m+1 \le i
\le 2m$.
\end{lemma}
{\em Proof.} We shall prove the assertion (1) by induction on $n
\in {\Zp}$.

For $n=1$ we can see directly  that $Q^{2} e^{\ga_i -\a}\ne 0$ (or see below).

Assume now that (1) holds for certain $n \in {\Zp}$. Since
$V_{L+\ga_i} \otimes F$ is a simple module for the simple vertex
operator superalgebra $V_L \otimes F$ we have that
$$ Y(E,z) Q^{2n} e^{\ga_i -n \a} \ne 0,$$
(for the proof see \cite{DL}).
So there is $j_0 \in {\Z}$ such that
$$E_{j_0} Q^{2n} e^{\ga_i - n\a} \ne 0 \quad \mbox{and} \quad
E_{j} Q^{2n} e^{\ga_i - n\a} = 0 \ \mbox{for} \ j > j_0.$$
Since
$$ E_{j_0} Q^{2n} e^{\ga_i - n\a} = \frac{1}{(n+1) (2n+1)} Q^{2n +2} (e^{-\a} _{j_0}
e^{\ga_i - n \a}), $$
we have that $j_0 \le i - 1 - (2m+1) n$. By using the fusion rules
from Proposition \ref{fusion}, we conclude that
$$ e^{-\a} _{j_0} e^{\ga_i - n \a} \in U ({\NS}). e^{\ga_i - (n+1)
\a}$$
 and therefore
$Q^{2n+2} e^{\ga_i - (n+1) \a} \ne 0$, which proves (1).
Notice that the idea used in the induction step, and fusion rules from Proposition \ref{fusion} can
be alternatively used to show that $Q^2 e^{\gamma_i-\alpha} \neq 0$.

The proof of (2) is similar so we omit it here.

\qed

\begin{remark} {\em It would be desirable - in parallel with the Virasoro algebra case - to have a direct proof of Lemma \ref{non-triv} with no reference to
fusion rules. However, the Virasoro algebra approach based on matrix coefficients does not apply verbatim
to superconformal $(1,2m+1)$-models, so we decided to give a proof which uses
the theory of vertex algebras and fusion rules. We found this approach to be quite elegant. We also remark that Iohara and Koga proved certain properties of screening operators among super Feigin-Fuchs modules in Theorem 3.1, \cite{IK} (see also \cite{MR}), but it is not clear
whether these results can be used to prove Lemma \ref{non-triv}.}
\end{remark}

As in the Virasoro algebra case the $N=1$ Feigin-Fuchs modules are
classified according to their embedding structure. For our purposes  we shall focus only on modules of certain type (Type 4
and 5 in \cite{IK}). These modules are either semisimple (Type 5)
or they become semisimple after quotienting with the maximal
semisimple submodule (Type 4). As usual the singular vectors will
be denoted by $\bullet$ and cosingular vectors with $\circ$.

The following result follows directly from Lemma \ref{non-triv}
and the structure theory of super Feigin-Fuchs modules  \cite{IK}
after some minor adjustments of parameters (cf. Type 4 embedding
structure).
\begin{theorem} \label{str-fock-ns-1}
Assume that $i \in \{0, \dots, m-1\}$.
 \item[(i)] As a  module for the Neveu-Schwarz algebra, $V_{L+\ga_i} \otimes F$  is generated by the family of
singular and cosingular vectors $ \widetilde{Sing}_{i}  \bigcup
\widetilde{CSing}_{i}$, where
$$  \widetilde{Sing}_{i} =  \{ u_i ^{(j, n)} \ \vert \ j, n \in {\N}, \ 0 \le j \le 2n \}; \
\
  \widetilde{CSing}_{i} =  \{ w_i ^ {(j, n)} \ \vert \ n \in {\Zp}, 0 \le j \le 2n-1 \}. $$
These vectors satisfy the following relations:
\bea
&&  u_i ^{(j, n)}= Q ^{j} e ^{\ga_i -n \alpha} , \ \ Q ^{j} w_i
^{(j, n)} = e ^{\ga_i + n \a}. \nonumber \eea
 The submodule generated by singular vectors
$\widetilde{Sing}_i$, denoted by $S\Lambda(i+1)$, is isomorphic to
\footnote{ In this section notation $k \ L^{\NS} (c,h)$ means
$L^{\NS} (c,h) ^{\oplus k}$, $k \in {\N}$.}
$$
\bigoplus_{n =0 } ^{\infty} (2n +1) L^{\NS}(c_{2m+1,1}, h^{2i+1,
2n+1}).  $$

\item[(ii)] For the quotient module we have
$$S \Pi(m-i):=(V_{L+\ga_i} \otimes F) / S\Lambda(i+1) \cong
  \bigoplus_{n =1 } ^{\infty}
(2 n ) L^{\NS}(c_{2m+1,1}, h^{2i+1, -2n +1}). $$

\end{theorem}

The situation described in Theorem \ref{str-fock-ns-1} can be
depicted by the following diagram:

 \be \label{type4.1}
\xymatrix{ j=-2 &  j=-1 & j=0  &  j=1 & j=2 \\  & & \bullet  &  & \\
& & \circ \ar[u] \ar[d] & \circ
\ar[d] & \\
& \bullet & \bullet & \bullet  & \\
& \circ \ar[u]\ar[d] & \circ \ar[u] \ar[d] & \circ \ar[u] \ar[d] &
\circ \ar[d] \\ \bullet & \bullet & \bullet & \bullet & \bullet
\\ \ar@{.}[u] & \ar@{.}[u] & \ar@{.}[u] & \ar@{.}[u] & \ar@{.}[u] } \ee

Let $M'$ be the contragradient $V$-module, where $V$ is a vertex
operator superalgebra. Then we have an isomorphism of $M(1)
\otimes F$-modules.
$$(M(1) \otimes e^{\frac{j}{2m+1}\alpha+i \alpha}
\otimes F)' \cong (M(1) \otimes e^{\frac{2m-j}{2m+1}\alpha-i
\alpha} \otimes F).$$ By taking direct sums we obtain the
following isomorphism of $\NS$-modules \be \label{dual-iso}
(V_{L+\gamma_i} \otimes F)' \cong V_{L+\gamma_{2m-i}} \otimes F.
\ee Since the dual functor interchanges cosingular and singular
vectors, Theorem \ref{str-fock-ns-1} implies the next result
(alternatively, use Type 4 embedding structure in \cite{IK}):

\begin{theorem} \label{str-fock-ns-2}

Assume that $i \in \{0, \dots, m-1 \}$.
 \item[(i)] As a  module for the Neveu-Schwarz algebra, $V_{L+\ga_{2m-i}} \otimes F$  is generated by the family of
singular and cosingular vectors $ \widetilde{Sing}'_{i}  \bigcup
\widetilde{CSing}'_{i}$, where
$$  \widetilde{Sing}'_{i} = \{ u_i^ {'(j, n)} \ \vert \ n \in {\Zp}, 0 \le j \le 2n-1
\}; \ \widetilde{CSing}'_{i} = \{ w_i^{'(j, n)} \ \vert \ j, n \in
{\N}, \ 0 \le j \le 2n \}. $$
These vectors satisfy the following relations:
\bea
&&  u_i ^{'(j, n)}= Q ^{j} e ^{\ga_{2m-i} -n \alpha} , \ \ Q ^{j}
w_i ^{'(j, n)} = e ^{\ga_{2m-i} + n \a}. \nonumber \eea

The submodule generated by singular vectors $\widetilde{Sing}_i$
is
%
%
is isomorphic to

$$ S\Pi(m-i)\cong \bigoplus_{n =1 } ^{\infty} (2 n ) L^{\NS}(c_{2m+1,1}, h^{2i+1,
-2n +1}).$$

\item[(ii)] For the quotient module we have
$$S \Lambda(i+1)\cong (V_{L+\ga_i} \otimes F) / S \Pi(m-i) \cong
  \bigoplus_{n =0 } ^{\infty}
(2 n+1 )  L^{\NS}(c_{2m+1,1}, h^{2i+1, 2n+1}).
$$

\end{theorem}

The embedding diagram for $V_{L+\gamma_{2m-i}} \otimes F$,
$i=0,...,m-1$ is now

\be \label{type4.2}
\xymatrix{ j=-1& j=0 & j=1  &  j=2  \\  & & \circ \ar[d]  &   \\
& \bullet & \bullet & \\
& \circ \ar[u] \ar[d] & \circ \ar[u] \ar[d] & \circ \ar[d]  \\
 \bullet & \bullet & \bullet & \bullet \\ \ar@{.}[u] & \ar@{.}[u] & \ar@{.}[u] & \ar@{.}[u] }
\ee

Finally, (\ref{dual-iso}) imply that $V_{L+\gamma_m} \otimes F$ is
a self-dual $V_L \otimes F$-module. In view of that, it is not surprising
that $V_{L+\gamma_m} \otimes F$ is a semisimple $\NS$-module.
More precisely, we have the following result
(for the proof see embedding
structure in Type 5 case in \cite{IK})

\begin{theorem} \label{str-fock-ns-3} As a  module for the Neveu-Schwarz algebra  $V_{L + \ga_{m}}
\otimes F$ is completely reducible and  generated by the family of
singular vectors
$$  \widetilde{Sing}_{m} =  \{ u_{m} ^{(j, n)} : = Q ^{j} e ^{\ga_{m} -n \alpha}  \ \vert \ j, n \in {\N}, \ 0 \le j \le 2n \}; $$

 and it is isomorphic to
$$
S\Lambda(m+1):=
 V_{L + \ga_{m}} \otimes F \cong  \bigoplus_{n=0} ^{\infty} (2n +1)L^{\NS}(c_{2m+1,1}, h^{2m+1, 2n
+1}).$$

\end{theorem}

The embedding structure in the last case is a totally disconnected
diagram

\be
 \label{type5}
\xymatrix{ j=-2 & j=-1 & j=0  &  j=1 & j=2 \\  & & \bullet  &  & \\
& \bullet & \bullet & \bullet & \\
\bullet & \bullet & \bullet & \bullet  & \bullet  \\
: & : & : & : & : } \ee

\section{The vertex operator superalgebra $ \overline{SM(1)}$}  \label{def-voa-u}

Let us fix a positive integer $m$. We shall first present the structure of the vertex operator
superalgebra $SM(1)$ as a module for the Neveu-Schwarz algebra. The next result
follows directly from Theorem \ref{str-fock-ns-1}.

\begin{theorem} \label{str-fock}
For every $n \in {\N}$,  set  $$u_n: = u_{0} ^{(n,n)} = Q^n
e^{-n\a}, \quad  w_{n+1} := w_{0} ^{(n+1,n+1)}.$$
 \item[(i)] The vertex operator superalgebra $SM(1)$, as a module for the vertex
operator superalgebra $L^{\NS}(c_{2 m+1,1},0)$, is generated by
the family of singular and cosingular vectors $ \widetilde{Sing}
\bigcup \widetilde{CSing}$, where
$$  \widetilde{Sing} =  \{ u_{n}  \ \vert \ n \in {\N}\}; \
\
  \widetilde{CSing} =  \{ w_{n} \ \vert \ n \in {\Zp} \}. $$
Moreover,  $U({\NS}) u_{n} \cong L^{\NS}(c_{2 m +1,1},
h^{1,2n+1})$ .
 \item[(ii)] The submodule generated by
vectors $u_n, n \in {\N}$ is isomorphic to $$ [\mbox{Sing}] \cong
\bigoplus_{n =0 } ^{\infty}
 L^{\NS}(c_{2 m +1,1},  h^{1,2n+1}).$$

\item[(iii)]The quotient module  is isomorphic to
$$M(1)/  [\mbox{Sing}] \cong \bigoplus_{n =1 } ^{\infty}
 L^{\NS}(c_{2 m +1,1},  h^{1,-2n+1}). $$
\item[(iv)] $Q u_0= Q {\vak} =0$, and $Q u_n \ne 0$, $ Q w_n \ne
0$ for every $n \ge 1$.
\end{theorem}

Our Theorem \ref{str-fock} immediately gives the following result.

\begin{proposition} \label{pomoc2}
 We have
$$L^{\NS} (c_{2 m +1, 1}, 0) \cong   W_0 =\mbox{Ker}_{SM(1)} Q  \ . $$
\end{proposition}

Define the following vertex algebra
$$\overline{SM(1)} = \mbox{Ker}_{SM(1) } \widetilde{Q}. $$

Since $ \widetilde{Q}$ commutes with the action of the
Neveu-Schwarz algebra, we have
$$L^{\NS}(c_{2m+1,1},0)\cong W_0 \subset  \overline{SM(1)}. $$
This implies that $ \overline{SM(1)} $ is  a vertex operator
subalgebra of $SM(1)$ in the sense of \cite{FHL} (i.e.,  $
\overline{SM(1)} $ has the same Virasoro element as $SM(1)$).

The following theorem will describe the structure of the vertex
operator superalgebra $\overline{SM(1)}$ as a
$L^{\NS}(c_{2m+1,1},0)$--module.

\begin{theorem} \label{generatori}
The vertex operator superalgebra $\overline{SM(1)}$ is isomorphic
to $ [\mbox{Sing}]$ as a $L^{\NS}(c_{2m+1,1},0)$--module, i.e.,
$$ \overline{SM(1)} \cong \bigoplus_{n =0} ^{\infty} L^{\NS}(c_{2 m +1,1}, h^{1,2n+1}).$$
\end{theorem}
{\em Proof.} By Theorem \ref{str-fock} we know that the
$L^{\NS}(c_{2m+1,1},0)$--submodule generated by the set
$\widetilde{Sing} $ is completely reducible. So to prove the
assertion, it suffices to show that the operator $ \widetilde{Q}$
annihilates   vector $v \in \widetilde{Sing} \cup
\widetilde{CSing} $  if and only if $v \in \widetilde{Sing} $.
Let $v \in  \widetilde{Sing} $, then $v = Q^{n} e ^{ -n \a}$ for
certain $n \in {\N}$. Since by Lemma \ref{rel-tilde}
$\widetilde{Q} e ^{ -n\a}=0$, we have
$$ \widetilde{Q} v =\widetilde{Q} Q^{n} e ^{ -n \a} = Q^{n}
\widetilde{Q} e ^{ -n\a} = 0. $$

Let now $v \in \widetilde{CSing}$. Then there is $n \in {\Zp}$
such that $Q ^{n} v = e ^{ n \a}$. Assume that $ \widetilde{Q} v =
0$. Then we have that
$$0=  Q ^{n} \widetilde{Q} v = \widetilde{Q} Q ^{n} v =
\widetilde{Q} e ^{ n \a}, $$
contradicting Lemma \ref{rel-tilde} (iii). This proves the
theorem.
 \qed

Next  we shall prove that the vertex operator algebra $
\overline{SM(1)} $ is generated by only two generators.

\begin{theorem} \label{singlet-generatori}
\item[(1)] The vertex operator superalgebra  $ \overline{SM(1)} $
is generated by $\tau$ and $H$.

\item[(2)] The vertex operator superalgebra $\overline{SM(1)}$ is
strongly generated by the set $$ \{ \tau, \omega, H,
G(-\tfrac{1}{2}) H \}.$$
\end{theorem}
{\em Proof.}
Let $U$ be the vertex subalgebra of $ \overline{SM(1)}$ generated
by $ \tau$ and $H$. We need to prove that $U = \overline{SM(1)}$.
Let $W_n$ by the (irreducible) ${\NS}$--submodule of $
\overline{SM(1)}$ generated by vector $u_n$. Then $W_n \cong
L^{\NS}(c_{2 m +1,1}, h ^{1,2n+1})$.
Using Lemma \ref{non-triv}  we see that
$$ \mbox{Ker}_{ \overline{SM(1)} }\  Q ^{n} \cong \bigoplus_{i = 0
} ^{n -1 } W_i .$$
To prove (1) suffices to show that $u_n \in U$ for every $n \in {\N}$. We
shall prove this claim by induction. By definition we have that
$u_0, u_1(= H) \in U$. Assume that we have $k \in {\N}$ such that
$u_{n} \in U$ for $n \le k$. In other words, the inductive
assumption is
$ \oplus_{i = 0} ^{k} W_i \subset U.$

We shall now prove that $u_{k+1} \in U$. Set $j = -(2 m +1)k -1$.
By Lemma \ref{non-triv}  we have
$$ Q ^{2 k + 2} e ^{- (k+ 1) \a} =  Q ^{2 k + 2} \left(e ^{-\a}
_{j} e ^{- k \a} \right)  \ne 0. $$
Next we notice that
$$ Q ^{k + 1}  ( H _j u_k ) = Q ^{k + 1} \left( Q e^{-\a}\right)
_j \left( Q^{k} e^{-k\a} \right)= \frac{1}{2 k + 1}  Q ^{2 k + 2}
\left(e ^{-\a} _{j} e ^{- k \a} \right), $$
which implies that
 $$Q ^{k + 1}  (H _j u_k) \ne 0.$$
  So we have found vector $H_j u_k
 \in U$ such that
 $$ \mbox{wt} ( H_j u_k)  = \mbox{wt} (u_{k+1}). $$
 This implies
 $$ H_j u_k \in \bigoplus_{i=0} ^{k+1} W_i \ \ \mbox{and} \ \ \
  H_j u_k \notin \bigoplus_{i=0} ^{k} W_i \ .$$
 Since $Q ^{k+1 } \left(\oplus_{i=0} ^{k} W_i\right)= 0$ and
 $ \mbox{wt} ( H_j u_k) = \mbox{wt} (u_{k+1})$  we conclude that
 there is a   constant $C$, $C \ne 0$,  such that
 $$ H_j u _k = C u_{k + 1} + u', \ \ u' \in \bigoplus_{i=0} ^{k}
 W_i \subset U. $$
 Since $H_j u_k \in U$, we conclude that $u_{k+1} \in U$.

 Therefore, the claim is verified, and the proof of (1) is complete.

 The proof of (1) shows  that $\overline{SM(1)}$ is spanned
by the vectors
\bea && u^{1}_{n_1}    \cdots  u^{r}_{n_r} {\bf 1}, \quad u^{i} \in
\{ \tau,   H \},  \label{spanning-set-singlet} \eea
such that for $1 \le i \le r$:
\bea && n _i \le -1 \quad \mbox{if} \ u^{i} =H \qquad \mbox{and}
\qquad n _i \le 0 \quad \mbox{if} \ u ^{i} = \tau.
\label{spanning-cond-singlet} \eea
 This implies that $\overline{SM(1)}$ is strongly generated by the set  $\{ \tau, \omega, H,
G(-\tfrac{1}{2}) H \}$, and (2) holds.  \qed

\ \vskip 5mm
 The following lemma   imply that for $i \ge -(2m +1)$
vectors $H_i H $ and $\widehat{H}_{i+1} \widehat{H} $ can be
constructed using only the action of the Neveu-Schwarz  operators
$L(n)$ and $G(n+\tfrac{1}{2})$ on the vacuum vector ${\vak}$.

\begin{lemma} \label{sing-z} We have:
 \bea && H_i H  \in W_0\cong L ^{\NS}( c_{2 m +1,1}, 0) \ \ \mbox{ for every} \  i \ge - (2 m
 +1),
 \nonumber \\
&&\widehat{H}_i \widehat{H}  \in W_0\cong L ^{\NS}( c_{2 m +1,1},
0) \ \ \mbox{ for every} \  i \ge - 2 m. \nonumber
 \eea
\end{lemma}

\begin{remark} {\em If we adopt notation used by physicists, then Theorem
\ref{singlet-generatori} implies that $\overline{SM(1)}$ is a
$\mathcal{W}(\frac{3}{2},2m+\frac{1}{2})$ superalgebra, meaning that
it is generated by primary fields of weight $\frac{3}{2}$ and
$2m+\frac{1}{2}$. In some physics papers
$\mathcal{W}(\frac{3}{2},2m+\frac{1}{2})$ super algebras are studied
by using general principles (e.g., Jacobi identities) but only for
low $m$. Because $\overline{SM(1)}$ shares many similarities with
the singlet algebra $\overline{M(1)}$ \cite{AdM-2007} we call
$\overline{SM(1)}$ {\em super singlet} vertex algebra.}
\end{remark}

\section{Zhu's algebra $A(\overline{SM(1)})$ and classification of
irreducible $\overline{SM(1)}$--modules}

In this section we completely determine Zhu's algebra
$A(\overline{SM(1)})$ and classify all irreducible
$\overline{SM(1)}$--modules. It turns out that the structure of
Zhu's algebra $A(\overline{SM(1)})$ is similar to the structure of
Zhu's algebra for $A(\overline{M(1)})$ studied in \cite{A-2003}
and  the proofs of the main results are completely analogous.

  Recall that $\widehat{H} = Q (e^{-\a}
\otimes \phi(-\tfrac{1}{2}))$. Clearly, $\widehat{H}$ is
proportional to $G(-\tfrac{1}{2}) H$ and therefore $\widehat{H}
\in \overline{SM(1)}$.

Next result shows that  Zhu's algebra $A(\overline{SM(1))}$ is
commutative.

\begin{theorem} \label{str-singlet-zhu}
Zhu's algebra $A(\overline{SM(1)})$ is spanned by the set
$$ \{ [\omega] ^{*s} [ \widehat{H}] ^{*t} \ \vert \ s,t \ge 0
\}.$$
In particular, Zhu's algebra $A(\overline{SM(1)})$ is isomorphic
to a certain quotient of the polynomial algebra ${\C}[x,y]$, where
$x$ and $y$ correspond to $[\omega]$ and $[\widehat{H}]$.
\end{theorem}

{\em Proof.} The proof follows from  Proposition \ref{abe},
Theorem \ref{singlet-generatori} and because $\tau$ and $H$
are odd vectors. \qed

Let $h^{r,s} = h_{2m+1,1} ^{r,s}$, so that $h^{2i+1,1}= \frac{ i
(i-2m)}{2 (2m+1)}$.

As in \cite{A-2003}, for $X=F$, $E$ or $H$ we let
$\widehat{X}(n):=\widehat{X}_{2m+n}$ (here as usual
$Y(\widehat{X},z)=\sum_{n \in \mathbb{Z}} {\widehat{X}}_n
z^{-n-1}$). In particular, $\widehat{H}(0)$ is a
  degree zero operator
acting on $\overline{SM(1)}$.   Since $\overline{SM(1)} \subset
M(1)$ every $M(1,\lambda) \otimes F$ is naturally an
$\overline{SM(1)}$-module.

Let ${\T}$ be the subspace of $M(1) \otimes F$ linearly spanned by
the vectors
 $$ a \otimes b, \quad \mbox{where} \ a \in M(1), \ b \in F, \
 \deg(b) > 0.$$
(So we only assume that $b$ is homogeneous in $F$ and that it is
not proportional to ${\bf 1}$).

The proof of the following  lemma is a consequence of the
definition of vertex superalegbra structure on $M(1) \otimes F$.

\begin{lemma}
 Let $\lambda \in \goth{h}^*$
 and $v_{\lambda}$ be the highest weight
vector in $M(1,\lambda) \otimes F$. Assume that $w \in {\T}$. Then
$o(w) v_{\l} = 0$.
\end{lemma}

 We have the following proposition
about the action of "Cartan subalgebra" of $\overline{SM(1)}$ on
the top component.

\begin{proposition} Let $\lambda \in \goth{h}^*$, $t=\langle
\alpha, \lambda \rangle$ and $v_{\lambda}$ the highest weight
vector in $M(1,\lambda) \otimes F$. Then we have
$$L(0) \cdot v_{\lambda}=\frac{t(t-2m)}{2(2m+1)} v_{\lambda},$$
$$\widehat{H}(0) \cdot v_{\lambda}= { t \choose 2m+1 } v_{\lambda}.$$
\end{proposition}
\noindent {\em Proof.} From the very definition of $Q$ and $H$ we
see that
$$\widehat{H}=\phi(1/2)S_{2m+1}(\alpha)\phi(-1/2) + w =S_{2m+1}(\alpha) + w,$$
where
$$w = S_{2m -1} (\a) \otimes \phi(-\tfrac{3}{2}) \phi(-\tfrac{1}{2}) + \cdots + {\bf 1} \otimes \phi(-2m -\tfrac{1}{2}) \phi(-\tfrac{1}{2})  \in {\T}.$$
 On the other hand it is known (cf. Proposition
3.1 in \cite{A}) that
$$S_r(\alpha)(0) v_{\lambda}={t \choose r}v_{\lambda}, \ \ r \geq 1.$$
The proof follows. $\Box$

It is not hard to see that $x(t)=\frac{t(t-2m)}{2(2m+1)}$ and
$y(t)={ t \choose 2m+1}$ parametrize the genus zero curve
$P(x,y)=0$ where \be \label{pxy} P(x,y)=y^2-C_m
\left(x+\frac{m^2}{2(2m+1)} \right) \prod_{i=0}^{m-1}
\left(x-\frac{i(i-2m)}{2(2m+1)} \right)^2, \ee where
$C_m=\frac{2^{2m+1}(2m+1)^{2m+1}}{(2m+1)!}$. Alternatively, notice
that we can write \be P(x,y)=y^2-C_m \prod_{i=0}^{2m}
\left(x-h^{2i+1,1} \right), \ee

By using arguments analogous to those in the proof of Lemma 6.1 from \cite{A-2003},
we obtain the following result:
\begin{lemma} \label{rel1-zhu}
In Zhu's algebra $A(\overline{SM(1)})$ we have the following
relation
$$[ \widehat{H}] * [ \widehat{H}]  = C_m \prod_{i=0} ^{2m}
([\omega]-h^{2i+1,1}),$$
where  $C_m$ is as above.
\end{lemma}

By using Theorem \ref{str-singlet-zhu},  Lemma \ref{rel1-zhu} and
the same proof as that of Theorem 6.1. from \cite{A-2003} we get:
\begin{theorem} \label{Zhu-singlet}
Zhu's algebra $A(\overline{SM(1)})$ is isomorphic to the
commutative, associative algebra ${\C}[x,y] / \la P(x,y) \ra$
where $\la P(x,y) \ra $ is the ideal in ${\C}[x,y]$ generated by
the polynomial
$$ P(x,y) = y^2 - C_m \prod_{i=0} ^{2m} (x- h^{2i+1,1}). $$
\end{theorem}

The fact that Zhu's  algebra $A(\overline{SM(1)})$ is commutative,
enable us to study irreducible lowest  weight representations of
the vertex operator superalgebra $\overline{SM(1)}$. For given
$(r, s) \in {\C} ^2$ such that $P(r,s) =0$ let ${\mathcal L}
({r,s})$ be the irreducible lowest weight
$\overline{SM(1)}$--module generated by the vector $v_{r,s} $ such that
$$ L(m) v = r\delta_{m,0}  v_{r,s} , \quad \widehat{H}(m) v =  s\delta_{m,0} r  v_{r,s} \quad  (m \ge 0).$$

Our Theorem \ref{Zhu-singlet} and standard Zhu's theory imply the following
classification result.

\begin{theorem} \label{clas-singlet}
The set
$$ \{ \mathcal{L}(r,s) \ \vert \ P(r,s) = 0 \}$$
provides all non-isomorphic irreducible
$\tfrac{1}{2}{\N}$-gradable $\overline{SM(1)}$-modules.
\end{theorem}

By using classification of irreducible $\overline{SM(1)}$--modules
and the same proof to that of Theorem 4.3 of \cite{AdM-2007} we
get:

\begin{corollary}
The vertex operator superalgebra $\overline{SM(1)}$ is simple.
\end{corollary}

\subsection{Logarithmic $\overline{SM(1)}$-modules}

In \cite{AdM-2007} we studied logarithmic modules for the singlet
vertex algebra $\overline{M(1)_p}$. Here we have a similar result.

As in \cite{AdM-2007}, let $M(1,\l) \otimes \Omega$ be an
$\hat{\goth{h}}$-module, where $\Omega$ is a two-dimensional vector
space and where $\alpha(0)|_{\Omega}$ is given by formula
\be \label{jordan} \left[\begin{array}{cc} \la \a , \lambda \ra & 1
\\ 0 & \la \a , \lambda \ra
\end{array}\right]
\ee in some basis $\{w_1,w_2 \}$ of $\Omega$
%
 (see also  \cite{M1}). Then
$M(1,\l) \otimes F \otimes \Omega$ carries an $\NS$-module
structure.

\begin{proposition} The vector space $M(1,\l) \otimes F \otimes
\Omega$, $\l \neq \frac{m}{2m+1} \alpha$ is a genuine logarithmic
$\overline{SM(1)}$-module \footnote{In other words, the module
involves nontrivial Jordan blocks with respect to the action of
$L(0)$.}, while for $\l=\frac{m}{2m+1} \alpha$, $M(1,\l) \otimes F
\otimes \Omega$ is an ordinary $\overline{SM(1)}$-module.
\end{proposition}

Notice that the previous result is in agreement with Theorem
\ref{Zhu-singlet}. More precisely, because of the linear term
$(x+\frac{m^2}{2(2m+1)})$ in $P(x,y)$, as in the proof of
Proposition 7.1 \cite{AdM-2007}, Theorem \ref{Zhu-singlet} can be
now used to show that there are no logarithmic self-extension of
$M(1,\frac{m}{2m+1}\alpha) \otimes F$.

\subsection{Further properties of $A(\overline{SM(1)})$}
In the next sections we shall make use of the following important
technical results.
\begin{proposition} \label{rel-zhu-singlet}
In Zhu's algebra $A(\overline{SM(1)})$ we have
$$ [Q^2 e^{-2\a}] = B_m    f_m ([\omega]) $$
where $$f_m([\omega]) =  \prod_{i=0} ^{3m} ([\omega] - h^{2i+1,1})
$$ and
$$B_m=  (-1) ^{m} \frac{ {2m \choose m} (2 (2m+1)) ^{3m+1} }{{4m
+1 \choose m } (3m+1) !
 ^{2}}.$$
\end{proposition}

{\em Proof.} First we notice that
$$ Q^{2} e^{-2\a} = {\nu} H_{-2m-2} H + v$$
where $ {\nu} \ne 0$ and $v \in U({\NS}).{\bf 1}$ (see also
\cite{AdM-triplet}, Lemma 3.3).
The above results on the structure of $A(\overline{SM(1)})$ implies
that
$$[Q^2 e^{-2\a}] = {\Phi}_m([\omega])
 $$
for certain ${\Phi}_m \in {\C}[x]$, $\deg {\Phi}_m \le 3 m+1$.
We shall evaluate the action of $Q^{2} e^{-2 \a}$ on top levels of
$\overline{SM(1)}$--modules $M(1,\l)\otimes F$. Let $v_{\l}$ be the
highest weight vector in
 $M(1,\l) \otimes F$.
First we notice that
$$ Q^{2} e^{-2\a} = \sum_{i=0} ^{4m +1} e^{\a}_{-i-1} e^{\a} _i
e^{-2\a}+ w, \quad \mbox{where} \ w \in {\T}. $$
 By using a direct calculation similar to that of \cite{AdM-triplet}  we see that
 \bea && o( Q^{2} e^{-2\a}) v_{\l} =
 \sum_{i=0} ^{\infty} o(e^{\a}_{-i-1}
e^{\a} _i e^{-2\a}) v_{\l} \nonumber \\
&& =\mbox{Res}_{z_1 } \mbox{Res}_{z_2} \sum_{i=0} ^{\infty}
z_1^{-i-1} z_2^i (z_1-z_2)^{2m+1} (z_1 z_2) ^{-4m-2} (1+z_1)^t
(1+z_2)^t v_{\l}  \nonumber \\
&& =\mbox{Res}_{z_1 } \mbox{Res}_{z_2}  (z_1-z_2)^{2m} (z_1 z_2)
^{-4m-2} (1+z_1)^t
(1+z_2)^t v_{\l}  \nonumber \\
 && = {\Phi}_m (\tfrac{1}{2 (2m+1)} (t ^{2} - 2 t m ) )v_{\l}= \widetilde{{\Phi}_m}(t)v_{\l}, \nonumber \\
 && \mbox{where} \quad \ \widetilde{\Phi_m }(t) =\sum_{k=0}
^{2 m } (-1) ^k  { 2 m \choose  k } { t \choose  4 m + 1 - k} {t
\choose
  2 m + 1+k} , \nonumber \\
 &&  t= \la \l, \a \ra . \nonumber \eea
 As in the proof of Lemma 3.4 in \cite{AdM-triplet} one can
prove the following identity
\bea \label{identitet} && \widetilde{\Phi_m }(t) =  {\bar A}_m { t
\choose 3 m +1} { t + m \choose 3m +1}, \quad \mbox{where}  \ {\bar
A}_m = \frac{(-1) ^{m} {2m \choose m} } {  {4m+1 \choose m} }. \eea
This implies
 $$  {\Phi}_m (\tfrac{1}{2 (2m+1)} (t ^{2} - 2  m t  )) =\widetilde{{\Phi}_ m }(t) = B_m f_m ( \tfrac{1}{ 2(2m+1)} (t ^{2} - 2 m t ) ).$$
Consequently, ${\Phi}_m$ is a non-trivial polynomial of degree
$3m+1$ and in $A(\overline{SM(1)})$ we have
 \bea \label{rel-q2}  && [Q^{2} e^{-2\a}] ={\Phi}_m ([\omega]) =
 B_m
 f_m ([\omega]), \ \ B_m \neq 0 \qed \eea

\vskip 5mm

 Define the following non-trivial vector $$U^{F,E}:=\mbox{Res}_z Y(F,z) E \frac{(z+1)^{2m }}{z} \in  \overline{SM(1)}. $$
 Set $U^{F,E} (0) := o(U^{F,E})= \sum_{i \ge 0}  {2m  \choose i} o( F_{i-1} E)$.

\begin{proposition} \label{rel-zhu-druga}
In Zhu's
algebra $A(\overline{SM(1)})$ we have:
$$ [U^{F,E}] = g([\omega]) [\widehat{H}]$$
where $g(x) \in \mathbb{C}[x]$ is of degree at most $m$.
\end{proposition}
{\em Proof.}
First we notice   that
$$ U^{F,E}= a H -  H\circ H $$
for certain $a \in U({\NS})$.
This implies that in Zhu's algebra $A(\overline{SM(1)})$, we have
\bea && [U^{F,E}] = g([\omega]) [\widehat{H}] \label{rel-pol}\eea
where $g\in {\C}[x]$ is a polynomial of degree at most $m$. (Here
we used the relation $[H\circ H]=0$, which holds in $A(\overline{SM(1)})$.)

\qed

It is not at all clear that $g(x)$ is a nonzero polynomial.

\section{The $N=1$ triplet vertex algebra $\mathcal{SW}(m)$ }

Define the following vertex superalgebra
$$\mathcal{SW}(m) = \mbox{Ker}_{V_L \otimes F} \widetilde{Q}.$$

Recall definition (\ref{def-hat}).
For any  $X \in \{E,F,H\}$,  $\widehat{X}$ is proportional to
$G(-\tfrac{1}{2}) X$, and therefore $\widehat{X} \in
\mathcal{SW}(m)$.

\begin{theorem} \label{generatori-triplet}
 \item[(1)]For every $m \ge 1$, $\mathcal{SW}(m)$ is an $N=1$
vertex operator superalgebra and $\mathcal{SW}(m) \cong
S\Lambda(1)$.

\item[(2)] The vertex operator superalgebra $\mathcal{SW}(m)$ is
generated by
$E$, $F$, $H$ and $\tau$.
\item[(3)] The vertex operator superalgebra $\mathcal{SW}(m)$ is
strongly generated by the set
$$\{\tau, \omega, E, F, H, \widehat{E}, \widehat{F}, \widehat{H}
\} . $$
\end{theorem}
\noindent {\em Proof.}
Recall the structure of $V_L \otimes F$ as a module for the
Neveu-Schwarz  algebra from Theorem \ref{str-fock-ns-1}. By using
Lemma \ref{rel-tilde} , similarly to the proof of Theorem
\ref{generatori}, we conclude  that
$\mathcal{SW}(m)$ is a completely reducible  module for the
Neveu-Schwarz algebra, generated  by the family of
singular vectors:
\bea \label{sing-trip} && Q^{j} e^{-n\a}, \quad n \in {\N}, \ j
\in \{0, \cdots, 2n \}. \eea
This proves (1).

  Let $Z_n$ be the Neveu-Schwarz  module generated by singular vectors
$$Q^{j} e^{-{\ell} \a}, \quad {\ell} \le n, \ j \in {\N} .$$
Therefore $\mathcal{SW}(m) = \bigcup_{n\in {\N} } Z_n$.
Let now ${\mathcal U}$ be the vertex subalgebra   of
$\mathcal{SW}(m)$ generated by $\tau, E, F, H$.
 Clearly,  ${\mathcal U} \subseteq \mathcal{SW}(m) $.
We shall prove that in fact  ${\mathcal U}=\mathcal{SW}(m) $.
In order to do so it is sufficient to show that $Z_n \subseteq {\mathcal U}$ for every $n \in {\Zp}$. We
shall prove this claim by induction on $n$. By the definition, the
claim holds for $n=1$. Assume now that $Z_n \subseteq {\mathcal
U}.$ Set $j_0 = (2m+1)n +1$. As in the proof of Theorem
\ref{singlet-generatori} we have
\bea
 F_{-j_0}  e^{-n \a} &=& e^{-(n+1) \a}, \nonumber \\
 E_{-j_0} Q^{2n}
e^{-n\a} &=& B_{2n+1} Q^{2n +2} e^{-(n+1)\a}, \nonumber \eea
{where} $B_{2n+1} \ne 0$ and
$$  H_{-j_0} Q^{j} e^{-n\a} = B_j  Q^{j+1} e^{-(n+1)\a} + v_j ' ,
$$
where $ v_j' \in Z_n, \ B_j \neq 0, \ 0 \le j \le 2n.$ These
relations imply that $Z_{n+1} \subseteq {\mathcal U}$. By induction
we conclude that $Z_n \subseteq {\mathcal U}$ for every $n \in
{\Zp}$ and therefore ${\mathcal U} = \mathcal{SW}(m) $. This proves
(2). The proof of (2) actually gives   that $\mathcal{SW}(m)$ is
spanned by the vectors
\bea && u^{1}_{n_1}    \cdots  u^{r}_{n_r} {\bf 1}, \quad u^{i} \in
\{ \tau, E, F, H \} \label{spanning-set} \eea
such that for $1 \le i \le r$:
\bea && n _i \le -1 \quad \mbox{if} \ u^{i} \in \{  E, F, H \}
\qquad \mbox{and} \qquad  n _i \le 0 \quad \mbox{if} \ u ^{i} =
\tau. \label{spanning-cond} \eea
 The assertion (3) follows.  \qed

\begin{theorem} Assume that $m \ge 1$. Then we have
\item[(1)] The vertex operator superalgebra $\mathcal{SW}(m)$ is
$C_2$--cofinite.
\item[(2)] The vertex operator superalgebra $\mathcal{SW}(m)$ is
irrational.
\end{theorem}
{\em Proof.} By using Proposition \ref{abe}, relation
(\ref{super-commute}) and Theorem \ref{generatori} we conclude that
$\mathcal{SW}(m) / C_2(\mathcal{SW}(m))$ is generated by
\bea \label{gen-c2}
\overline{\tau}, \overline{\omega}, \overline{E},
\overline{\widehat{E}}, \overline{F}, \overline{\widehat{F}},
\overline{H}, \overline{\widehat{H}},
 \eea
and that every two generators either commute or anti-commute. In
order to prove $C_2$--cofiniteness it suffices to prove that every
generator (\ref{gen-c2}) is nilpotent in $\mathcal{SW}(m) /
C_2(\mathcal{SW}(m))$.
Let $X$ be either $E$ or $F $. From Lemma \ref{pomoc1} we see that
$X_{-1}X=0$, and thus $\overline{X} ^{2}=0$. By using
$$G(-i-1/2)^2=L(-2i-1) \in U(\NS)$$
we get  $\overline{\tau} ^2 = 0$.  Similarly, from
$$H_{-1}H \in U(\NS) \cdot {\bf 1},$$
$$H_{-1}H=\sum_{
   k/ 2 + i_1+\cdots+i_k  + j_1 + \cdots + j_s = 4m +1} a_{i_1,...,i_k} G(-i_1-1/2) \cdots G(-i_k-1/2) L(-j_1 ) \cdots
L(-j_s){\bf 1},$$ { where} $$i_1 > i_2 > \cdots >i_k \geq 1, \
j_1, \dots, j_s \ge 2, \  a_{i_1,...,i_k} \in \mathbb{C},$$  it
follows that $H_{-1}H \in C_2(\striplet)$, and thus
$$\overline{H}^2=\overline{H_{-1}H}=0.$$

We also have $X_{-1} X = 0 $, $X \in \{\widehat{F},\widehat{E} \}$
in $\mathcal{SW}(m)$ (cf. Lemma \ref{nilp-hat}), so that
$$ \overline{X} ^{2} = 0 \ \ \mbox{in} \ \ \mathcal{SW}(m) / C_2 (
\mathcal{SW}(m) ). $$
Thus, it remains to prove that $\overline{\omega}$ and $\overline{
\widehat{H}}$ are nilpotent.  We prove this as in
\cite{AdM-triplet}. Since
$$ \widehat{E}_{-1} \widehat{F} + \widehat{F}_{-1} \widehat{E} + 2
\widehat{H}_{-1} \widehat{H} = 0$$
we get
$$ \overline{\widehat{H}} ^{4} = 0. $$
Moreover, the description of Zhu's algebra from Theorem
\ref{Zhu-singlet} implies that
$$ \overline{\widehat{H}} ^{2} = C_m \overline{\omega} ^{2m+1}, \quad (C_m \ne 0),$$
which implies that $\overline{\omega} ^{4m+2}=0$. Therefore, every
generator of $\mathcal{SW}(m) / C_2(\mathcal{SW}(m))$ is nilpotent
and    $\mathcal{SW}(m)  $ is $C_2$--cofinite. This proves (1).

Assertion (2) follows from the fact that $V_L \otimes F$ is not
completely reducible, viewed as $\mathcal{SW}(m)$--module.
 \qed

\section{Classification of irreducible $\mathcal{SW}(m)$--modules}

From the definition of Zhu's algebra and the structure of the
vertex operator superalgebra $\mathcal{SW}(m)$ follows:

\begin{proposition} \label{striplet-gen}
The associative algebra $A(\mathcal{SW}(m))$ is generated by
$[\widehat{E}]$, $[\widehat{H}]$, $[\widehat{F}]$ and $[\omega]$.
\end{proposition}
{\em Proof.} The proof follows from  Proposition \ref{abe},
Theorem \ref{generatori-triplet} and the fact that $\tau$, $E$,
$F$ and $H$ are all odd. \qed

\begin{theorem} \label{vazna-rel}
In Zhu's algebra $A(\mathcal{SW}(m))$ we have the following
relation
$$f_m([\omega]) = 0$$
where
$$f_m(x) = \prod_{i=0} ^{3m} (x-h^{2i+1,1}). $$
\end{theorem}
{\em Proof.} Since $O(\overline{SM(1)}) \subset
O(\mathcal{SW}(m))$, the embedding $\overline{SM(1)}  \subset
\striplet$ induces an algebra homomorphism $A(\overline{SM(1)})
\rightarrow A(\striplet)$. Applying this homomorphism to
Proposition \ref{rel-zhu-singlet}  and using the fact that $Q^{2}
e^{-2\alpha} \in O(\striplet)$  we get that $f_m([\omega]) = 0$ in
$A(\striplet)$.
 \qed

Alternatively, we can write the polynomial $f_m(x)$ as \be
\label{alter-fm} f_m(x) = (x-h^{2m+1,1})\prod_{i=0} ^{m-1}
(x-h^{2i+1,1})^2 \prod_{i=2m+1}^{3m}(x-h^{2i+1,1}), \ee indicating
possibility of existence of logarithmic modules of generalized
lowest conformal weight $h^{2i+1,1}$, $i=0,...,m-1$.

\begin{theorem} \label{ireducibilni-moduli}
\item[(1)] For every $0 \leq i \leq m$,  $S\Lambda(i+1)$ is an
irreducible $\tfrac{1}{2}\N$--gradable $\mathcal{SW}(m)$--module, with
the top component $S\Lambda(i+1) (0)$ of  lowest weight
$h^{2i+1,1}$. Moreover, $S\Lambda(i+1) (0)$ is an $1$--dimensional
irreducible $A(\mathcal{SW}(m))$--module.

\item[(2)]  For every $0 \leq j \leq m-1$ $,  S\Pi(m -j)$ is an
irreducible $\tfrac{1}{2}\N$--gradable $\mathcal{SW}(m)$--module, with
the top component $S\Pi(m-j) (0)$ of lowest weight $h^{2 i +1,1}$
where $i=2m+1 +j$. Moreover, $S\Pi(m-j) (0)$ is an
$2$--dimensional irreducible $A(\mathcal{SW}(m))$--module.
\end{theorem}
{\em Proof.} Proof is similar to that of Theorem 3.7 in
\cite{AdM-triplet} so we omit it here. \qed

Applying the previous theorem in the case of $\mathcal{SW}(m)
=S\Lambda(1)$ we get:
\begin{corollary}
The vertex operator superalgebra $\mathcal{SW}(m)$ is simple.
\end{corollary}

As in \cite{AdM-triplet} we have the following result
\begin{proposition} \label{zhu-komut}
In Zhu's associative algebra we have \bea &&
[\widehat{H}]*[\widehat{F}]-[\widehat{F}]*[\widehat{H}]=-2q([\omega])[\widehat{F}], \\
&& [\widehat{H}]*[\widehat{E}]-[\widehat{E}]*[\widehat{H}]=2q([\omega])[\widehat{E}] \\
&&
[\widehat{E}]*[\widehat{F}]-[\widehat{F}]*[\widehat{E}]=-2q([\omega])[\widehat{H}].
\eea where $q$ is a certain polynomial. \end{proposition}

\begin{theorem} \label{class-ired} The set $$\{ S\Pi(i)(0) : 1 \leq i \leq m  \} \cup
\ S\Lambda(i)(0): 1 \leq i \leq m+1 \}$$ provides, up to
isomorphism, all irreducible modules for  Zhu's algebra
$A(\mathcal{SW}(m))$.
\end{theorem}
{\em Proof.}
The proof  is similar to that of Theorem 3.11 in
\cite{AdM-triplet}.
 Assume that $U$ is an irreducible $A(\mathcal{SW}(m))$--module. Relation
$f_m ([\omega]) = 0$ in $A(\mathcal{SW}(m))$ implies that
$$L(0) \vert U = h^{2i+1,1} \ \mbox{Id}, \quad \mbox{for} \quad i \in \{0, \dots, m \} \cup \{2m+1, \dots, 3m\}.$$

Assume first that $i= 2m+1 + j$ for $ 0 \le j \le m-1$. By
combining Propositions \ref{zhu-komut} and Theorem
\ref{ireducibilni-moduli} we have that $q(h^{2i+1,1}) \ne 0$.
Define
$$ e= \frac{ 1}{ \sqrt{2} q(h^{2 i+1,1})} [\widehat{E}], \quad f= -\frac{1}{ \sqrt{2}q(h^{2i+1,1})} [\widehat{F}], \quad h= \frac{1}{q(h^ {2i+1,1})} [\widehat{H}] .$$

Therefore $U$ carries the structure of an irreducible,
$\goth{sl}_2$--module with the property that $e^2 = f^2 = 0$  and
$h \ne 0$ on $U$. This easily implies that $U$ is a $2$--dimensional
irreducible $\goth{sl}_2$--module. Moreover,  as an
$A(\mathcal{SW}(m))$--module $U$ is isomorphic to  $S\Pi(m-j)(0)$.

Assume next that $ 0 \le i \le m$. If $q(h^{2i+1,1}) \ne 0$, as
above
 we conclude that $U$ is an irreducible $1$--dimensional $sl_2$--module. Therefore $U \cong S\Lambda(i+1)(0)$.

 If $q(h^{2i+1,1})=0$,
 from Proposition \ref{zhu-komut} we have that the action of generators of $A(\mathcal{SW}(m))$ commute on $U$. Irreducibility of
  $U$ implies that $U$ is $1$-dimensional. Since $[\widehat{H}], [\widehat{E}] ^2, [\widehat{F}] ^2$ must act trivially on $U$, we conclude that $[\widehat{H}],
   [\widehat{E}] , [\widehat{F}] $ also act trivially on $U$. Therefore
  $U \cong S\Lambda(i+1)(0)$.
  \qed

As a consequence of the previous theorem we have.

\begin{theorem} \label{class-ired-modules} The set $$\{ S\Pi(i) : 1 \leq i \leq m  \} \cup
\{ S\Lambda(i)  : 1 \leq i \leq m +1 \}$$ provides, up to
isomorphism, all irreducible modules for the vertex operator
superalgebra $\mathcal{SW}(m)$.
\end{theorem}

\section{On the structure of Zhu's algebra $A(\striplet)$}

As in \cite{AdM-triplet}, the main difficulty in description of Zhu's algebra $A(\striplet)$ is  that of not having
a good understanding of  logarithmic $\striplet$-modules. For
the triplet $\mathcal{W}(p)$ this problem can be resolved, at least if $p$ is prime, by
using modular invariance. We believe  the same approach can be
applied for $\mathcal{SW}(m)$, which would require a super version
of Miyamoto's result \cite{Miy}.  This is the main reason why in this part we
focus mostly on the case $2m+1$ is prime, but we expect all
results to be true in general.

 In many ways this section is analogous to Section 5
(and Appendix) in \cite{AdM-triplet}, but as we shall see there are some important differences.

First a few generalities regarding the Lagrange interpolation
polynomial.

\begin{proposition} \label{interp} Let $S=\{  (x_1,y_1),...,(x_n,y_n) \}$, $x_i \neq x_j$ be a set of points
in $\mathbb{C}^2$ such that their Lagrange interpolation
polynomial $L_n(x)$ is of degree exactly $n-1$. Then every
interpolation polynomial of degree exactly $n$ is given by
$$Q_\lambda(x)=L_n(x)+\lambda \prod_{i=1}^{n}(x-x_i), \ \ \lambda
\neq 0.$$
\end{proposition}
\noindent {\em Proof.} Let $P(x)$ be an arbitrary interpolation
polynomial of degree $n$.  Then for some $\lambda$, the polynomial
$P(x)-\lambda \prod_{i=1}^n(x-x_i)$ is of degree less or equal
$n-1$, but not zero. But then $P(x)-\lambda \prod_{i=1}^n
(x-x_i)=L_n(x).$ \qed

\begin{lemma} \label{interp-1}
Let $L_m(x)$ be the Lagrange interpolation polynomial for
$(h^{2i+1,1},{i \choose 2m+1})$, where $ 2m+1 \leq i \leq 3m$. If
we let $r(t)=L_n(\frac{t(t-2m)}{2(2m+1)})$, then
$$r(t)=\frac{{\prod_{i=2m+1}^{3m}}
(t-i)(t-2m+i)}{(2m+1)!}$$
$$\times \sum_{i=2m+1} ^{3m} \left(\frac{(i!)^2(-1)^{i+m}}{(i-2m-1)!^2(3m-i)!(i+m)!}(\frac{1}{t-i}-\frac{1}{t-2m+i}) \right) \in \mathbb{C}[t].$$
\end{lemma}

Now,  we have an important technical result (in a slightly
different setup a similar result has been proven in Appendix of
\cite{AdM-triplet}).

\begin{proposition} \label{interp-2} For every $m \geq 1$ we have
$$L_m(h^{2i+1,1}) \neq 0, \ \ 0 \leq i \leq m.$$
\end{proposition}
\noindent {\em Proof.} As in \cite{AdM-triplet} it suffices to let
$$s(t)=\frac{r(t)}{\prod_{i=2m+1}^{3m} (t-i)(t-2m+i)}$$
and check first
$$s(0)<0, \ \ s(1) <0,$$
which follows by using hypergeometric summations. That
$r(h^{2i+1,1}) \neq 0$ for $0 \leq i \leq m$ follows now from the
recursion
$$s(t)(m+t)(2m+1-t)^2=2(m+1-t)(2m^2+2tm-2-t^2+2t)s(t-1)+(t-1)^2(3m+2-t)s(t-2),$$
because all coefficients in the recursion are positive for $1 \leq
t \leq m$. \qed

As in Appendix of \cite{AdM-triplet} we now observe that
$$\widehat{H}*\widehat{F}=a. F,$$
where
$$a \in U(\goth{ns}).$$
From
$${\rm deg}(\widehat{H}_{-1} \widehat{F})=4m+2,$$
and
\be \label{hf}
[\widehat{H}]*[\widehat{F}]=-q([\omega])[\widehat{F}], \ee
 for
some $q \in \mathbb{C}[x]$. It follows that $q([\omega])$ is a
polynomial of degree at most $m$. In \cite{AdM-triplet} this
observation was sufficient to argue that $q$ has to be the
interpolation polynomial. However, in view of Proposition
\ref{interp} and Lemma \ref{interp-1}, we are unable to argue that
$q=L_m$, because $L_m$ is of degree $m-1$. Thus, it is not clear what the $q$
polynomial should be.

\begin{proposition} \label{rel-zhu-prop} Let $g(x)$ be as in Proposition \ref{rel-zhu-druga}
and  $$u(x)=\prod_{i=2m+1}^{3m} (x-h^{2i+1,1}).$$
Then
$$g(x)=D_m u(x),$$
for some constant $D_m$. Moreover,
\be \label{conj-u}  D_m u([\omega]) *[\widehat{X}]=0, \  \ X \in
\{F,H,E \}. \ee
\end{proposition}
{\em Proof.} First we notice that $U^{F,E}= F \circ E \in
O(\striplet)$. Then Proposition \ref{rel-zhu-druga} implies that
$$ g([\omega]) * [\widehat{H} ] = 0 \quad \mbox{in} \quad
A(\striplet)$$ for some polynomial of degree at most $m$. Because
we already know all irreducible $\mathcal{SW}(m)$-modules we also
know that $g([\omega])$ must act as zero on all
$\mathcal{SW}(m)$-modules with two-dimensional highest weight
subspaces (here $[\widehat{H}]$ acts nontrivially). Thus we know
that
$$g([\omega])=D_m u([\omega])$$
for some constant $D_m$. Since $Q$ preserves $O(\striplet)$ we get
(\ref{conj-u}). \qed

It is crucial for our considerations to show that $D_m \neq 0$
(i.e., $g(x) \neq 0$). This will requires an explicit computation
of $U^{F,E}(0)$ on the top degree subspaces of certain
$\overline{SM(1)}$-modules. We have the following result.

\begin{theorem} \label{2m1prime} If $m \in \mathbb{N}$ such that $2m+1$ is a prime integer, then  $g(x)  \neq 0$.
\end{theorem}

For the proof of this important technical result we refer the reader to Appendix.

If $D_m \neq 0$,  Proposition \ref{rel-zhu-prop} and
(\ref{hf}) we get
$$[\widehat{H}]*[\widehat{F}]=-q([\omega])[\widehat{F}]=-q'([\omega])[\widehat{F}],$$
where $q'([\omega])$ is a polynomial of degree $m-1$, which forces
$q'=L_m$. We should say here that in \cite{AdM-triplet} the
formula (\ref{conj-u}) was a consequence of a formula analogous to
(\ref{hf}).

\begin{theorem} \label{zhu-relations}
Assume that $2m+1$ is prime or $D_m \neq 0$. Then we have

\begin{itemize}
\item[(i)]   $[\widehat{E}]^2=[\widehat{F}]^2=0$

\item[(ii)] $[\widehat{H}]^2=C_m P([\omega])$, where
$$P(x)=\prod_{i=0}^{2m}(x-h^{2i+1,1}) \in \mathbb{C}[x]$$
and $C_m$ is a nonzero constant.

\item[(iii)]
$$[\widehat{H}]*[\widehat{F}]=-[\widehat{F}]*[\widehat{H}]=-q([\omega])*[\widehat{F}],$$
$$[\widehat{H}]*[\widehat{E}]=-[\widehat{E}]*[\widehat{H}]=q([\omega])*[\widehat{E}],$$ where $q(x)$ is a nonzero
polynomial of degree $m-1$ and $$q(h^{2i+1,1}) \neq 0, \ \ 0 \leq
i \leq m.$$

\item[(iv)] \bea &&
[\widehat{H}]*[\widehat{F}]-[\widehat{F}]*[\widehat{H}]=-2q([\omega])[\widehat{F}], \nonumber \\
&& [\widehat{H}]*[\widehat{E}]-[\widehat{E}]*[\widehat{H}]=2q([\omega])[\widehat{E}], \nonumber \\
&&
[\widehat{E}]*[\widehat{F}]-[\widehat{F}]*[\widehat{E}]=-2q([\omega])[\widehat{H}],
\nonumber \eea where $q(x)$ is as in (iii).

\item[(v)] $$ \prod_{i=2m+1}^{3m}([\omega]-h^{2i+1,1})*[X]=0, \ \
X \in \{ \widehat{E},\widehat{F}, \widehat{H} \}.$$

\item[(vi)] The center of $A(\striplet)$ is a subalgebra generated
by $[\omega]$.

\end{itemize}
\end{theorem}

\noindent {\em Proof.} We recall that $\striplet$ is generated by
$[\omega]$ and $[\widehat{X}]$, $X=F$, $H$ and $E$ (see
Proposition \ref{striplet-gen}).

For (i) we recall \cite{AdM-triplet} that $Q$ lifts to a
derivation of $A(\striplet)$, denoted by the same symbol. Now,
because of Lemma \ref{nilp-hat} we have
$$[\widehat{F}]*[\widehat{F}]=[\widehat{E}]*[\widehat{E}]=0.$$

Part (ii) has been proven in Lemma \ref{rel1-zhu}.

It is left to show relations (iii), (iv) and (v). As in
\cite{AdM-triplet} we compute
$$0=Q([\widehat{F}]*[\widehat{F}])=[\widehat{H}]*[\widehat{F}]+[\widehat{F}]*[\widehat{H}],$$
which yields
$$[\widehat{H}]*[\widehat{F}]=-[\widehat{F}]*[\widehat{H}].$$
After an application of $Q^2$ on the previous equation we get
$$[\widehat{H}]*[\widehat{E}]=-[\widehat{E}]*[\widehat{H}].$$

Two remaining formulas in (iii)  \be \label{hfq}
[\widehat{H}]*[\widehat{F}]=-q([\omega])*[\widehat{F}], \ee
$$[\widehat{H}]*[\widehat{E}]=q([\omega])*[\widehat{E}],$$
have already been proven in the discussion preceding the theorem.

The relation (iv) follow from (iii) (cf. \cite{AdM-triplet}). Part
(v) follows directly from Proposition \ref{rel-zhu-prop}. Part
(vi) follows from the fact that $q([\omega])$ is a unit in
$A(\striplet)$. \qed

\begin{corollary}
Under the assumptions of Theorem \ref{zhu-relations}, the
associative algebra $A(\striplet)$ is spanned by
$$\{ [\omega]^i, \ 0 \leq  i  \leq  3m \} \cup \{[\omega]^i*[X], 0 \leq i
\leq m-1, \ X=\widehat{E},\widehat{F} \ \ {\rm or} \ \widehat{H}
\}.$$ Thus, $A(\striplet)$ is at most $6m+1$-dimensional.
\end{corollary}

By using the same ideas as in \cite{AdM-triplet} it is not hard to
show that

\begin{theorem} Assume  $2m+1$ be prime or $D_m \neq 0$. Then the Zhu's algebra decomposes as a sum of ideals
$$A(\striplet)=\bigoplus_{i=2m+1}^{3m}  \mathbb{M}_{h^{2i+1,1}} \oplus \bigoplus_{i=0}^{m-1} \mathbb{I}_{h^{2i+1,1}}
\oplus \mathbb{C}_{h^{2m+1,1}},$$ where $\mathbb{M}_{h^{2 i+1,1}}
\cong M_2(\mathbb{C})$, $1 \leq {\rm dim}(\mathbb{I}_{h^{2
i+1,1}}) \leq 2$ and $\mathbb{C}_{h^{2m+1,1}}$ is one-dimensional.
\end{theorem}

It is also not hard to find explicit generators for every ideal,
in parallel with \cite{AdM-triplet}.

As with the triplet we expect that all $\mathbb{I}_{h^{2i+1,1}}$
are two-dimensional (which is related to existence of logarithmic
modules). This is equivalent to

\begin{conjecture} \label{conj-zhu-6m1}
The associative algebra $A(\striplet)$ is $6m+1$-dimensional. Then
the center of $A(\striplet)$ is $3m +1$-dimensional.
\end{conjecture}

\begin{remark} \label{dong} {\em Dong and Jiang have recently proven \cite{DJ} that if $A(V)$ is
semisimple and every irreducible admissible module is an ordinary
module, then $V$ is rational. It is feasible to assume that their
result applies for vertex operator superalgebras. This would imply
${\rm dim}(\mathbb{I}_{h^{2i+1,1}})=2$ for at least one $i$, and
in particular
$${\rm dim} \ A(\mathcal{SW}(1))=7.$$
(Note that in the case $m=1$, $D_1 \neq 0$ certainly holds.)}
\end{remark}

\section{Modular properties of characters of irreducible $\mathcal{SW}(m)$-modules}

We first introduce several  basic facts regarding classical modular forms needed for description
of irreducible $\striplet$ characters. The Dedekind
$\eta$-function is usually defined as the infinite product
$$\eta(\tau)=q^{1/24} \prod_{n=1}^\infty (1-q^n),$$
an automorphic form of weight $\frac{1}{2}$. As usual in all these
formulas $q=e^{2 \pi i \tau}$, $\tau \in \mathbb{H}$ \footnote{Here $\tau$ - the coordinate of $\mathbb{H}$ - should not be confused with the superconformal vector used in previous sections.}. We also
introduce \bea \label{weber} \goth{f}(\tau) &=& q^{-1/48}
\prod_{n=0
}^\infty (1+q^{n+1/2}), \\
\goth{f}_1(\tau) &=& q^{-1/48} \prod_{n=1}^\infty (1-q^{n-1/2}), \\
\goth{f}_2(\tau) &=& q^{1/24} \prod_{n=1}^\infty (1+q^n).
\eea These (slightly normalized) Weber functions form a vector-valued modular form of
weight zero. More precisely,
$$\goth{f}(-1/\tau)=\goth{f}(\tau) ,\  \goth{f}_2(-1/\tau)=\frac{1}{\sqrt{2}}\goth{f}_1(\tau), \  \goth{f}_1(-1/\tau)=\sqrt{2} \goth{f}_2(\tau), $$
$$\goth{f}(\tau+1)=e^{-2 \pi i/48} \goth{f}_1(\tau), \  \goth{f}_2(\tau+1)=e^{2 \pi i/24} \goth{f}_2(\tau), \  \goth{f}_1(\tau+1)=e^{-2 \pi i/48}  \goth{f}(\tau). $$

In what follows, we denote by $$\Theta_{j,k}(\tau)=\sum_{n \in
\mathbb{Z}} q^{(2kn+j)^2/4k}$$ Jacobi-Riemann $\Theta$-series
where $j \in \mathbb{Z}$ and $k \in \mathbb{N}/2$. We also let
$$(\partial \Theta)_{j,k}(\tau)=\sum_{n \in \mathbb{Z}}
(2kn+j)q^{(2kn+j)^2/4k}.$$ Then we have transformation formulas
(notice that here $k \in \mathbb{N}/2$ so $\Theta_{j,k}(\tau)$ is
not invariant under $\tau \longrightarrow \tau +1$ in general):
\bea \label{mod-transf-1} && \eta(-1/\tau)=\sqrt{-i\tau}
\eta(\tau), \ \ \eta(\tau+1)=e^{\pi i/12} \eta(\tau) \\
&& \Theta_{j,k}(-1/\tau)=\sqrt{\frac{-i
\tau}{2k}}\sum_{j'=0}^{2k-1} e^{i \pi j j'/k} \Theta_{j',k}(\tau) \\
&& \Theta_{j,k}(\tau+2)=e^{i \pi j^2/k}\Theta_{j,k}(\tau) \\
&& (\partial \Theta)_{j,k}(\tau+2)=e^{i \pi j^2/k}(\partial
\Theta)_{j,k}(\tau), \\
&& (\partial \Theta)_{j,k}(-1/\tau)=(- \tau)\sqrt{-i \tau/2k}
\sum_{j'=1}^{2k-1} e^{i \pi j j'/k} (\partial
\Theta)_{j',k}(\tau). \eea

For a vertex operator algebra module $M$ we define its
graded-dimension or simply character
$$\chi_{M}(\tau)={\rm tr}|_{M} q^{L(0)-c/24}.$$
If $V=L^{\NS}(c_{2m+1,0},0)$ and $M=L(c_{2m+1,0},h^{2i+1,2n+1})$,
then (see \cite{IK}, for instance)
 \be \label{irr-ns}
 \chi_{L^{\NS}(c_{2m+1,1},h^{2i+1,2n+1})}(\tau)=q^{\frac{m^2}{2(2m+1)}}\frac{\goth{f}(\tau)}{\eta(\tau)}\left(q^{h^{2i+1,2n+1}}-q^{h^{2i+1,-2n-1}}\right).
 \ee

By combining Theorem \ref{str-fock-ns-1}, \ref{str-fock-ns-2} and
\ref{str-fock-ns-3}, and formula (\ref{irr-ns}) we obtain

\begin{proposition} \label{char-irr-mod} For $i=0,...,m-1$
\bea \label{sl-char} && \chi_{S
\Lambda(i+1)}(\tau)=\frac{\goth{f}(\tau)}{\eta(\tau)}\left(\frac{2i+1}{2m+1}
\Theta_{m-i,\frac{2m+1}{2}}(\tau)+\frac{2}{2m+1}(\partial
\Theta)_{m-i,\frac{2m+1}{2}}(\tau) \right), \\
&& \label{pi-char} \chi_{S
\Pi(m-i)}(\tau)=\frac{\goth{f}(\tau)}{\eta(\tau)}\left(\frac{2m-2i}{2m+1}
\Theta_{m-i,\frac{2m+1}{2}}(\tau)-\frac{2}{2m+1}(\partial
\Theta)_{m-i,\frac{2m+1}{2}}(\tau)\right). \eea Also, \be \chi_{S
\Lambda(m+1)}(\tau)=\frac{\goth{f}(\tau)}{\eta(\tau)}
\Theta_{0,\frac{2m+1}{2}}(\tau). \ee
\end{proposition}

For purposes of modular invariance, it is also important to
compute {\em supercharacters} of irreducible modules. Let us
recall that a supercharacter of a $V$-module $M$ is defined
$$\chi_{M}^F(\tau)={\rm tr}|_{M} \sigma q^{L(0)-c/24},$$
where $\sigma$ is the sign operator taking values $1$ (resp. $-1$)
on even (resp. odd) vectors.

In parallel with Proposition \ref{char-irr-mod}, it is not hard to
compute irreducible supercharacters of $\striplet$-modules. Here
is an explicit description in terms of $\Theta$-constants and
their derivatives.

\begin{proposition} \label{superchar-irr-mod} For $i=0,...,m-1$
\bea \label{sl-schar} && \chi^F_{S \Lambda(i+1)}(\tau)= \nonumber \\
&& \frac{\goth{f}_2(\tau)}{\eta(\tau)} \biggl(
\frac{2i+1}{2m+1}\left( \Theta_{2(m-i),2(2m+1)}(
\tau)-\Theta_{2(m+i+1),2(2m+1)}( \tau)\right)+ \\
&& +\frac{1}{2m+1} \left( (\partial \Theta)_{2(m-i),2(2m+1)}(\tau)
-(\partial \Theta)_{2(m+i+1),2(2m+1)}(\tau)\right)\biggr). \eea

\bea \label{p1-schar} && \chi^F_{S \Pi(m-i)}(\tau)= \nonumber \\
&& \frac{\goth{f}_2(\tau)}{\eta(\tau)} \biggl(
\frac{2m-2i}{2m+1}\left( \Theta_{2(m-i),2(2m+1)}(
\tau)-\Theta_{2(m+i+1),2(2m+1)}( \tau)\right)+ \\
&& -\frac{1}{2m+1} \left( (\partial \Theta)_{2(m-i),2(2m+1)}(\tau)
-(\partial \Theta)_{2(m+i+1),2(2m+1)}(\tau)\right)\biggr). \eea

Also, \be \chi_{S
\Lambda(m+1)}^F(\tau)=\frac{\goth{f}_2(\tau)}{\eta(\tau)}
\biggl(\Theta_{0,2(2m+1)}(\tau)-\Theta_{2(2m+1),2(2m+1)}(\tau)
\biggr). \ee
\end{proposition}

As in \cite{F} we now study modular invariance properties of
irreducible $\striplet$ characters and supercharacters. We only
consider some special modular transformations. For example,
$$\chi_{S
\Lambda(i+1)}(-1/\tau)=
\frac{\goth{f}(\tau)}{\eta(\tau)}\sum_{k=0}^{2m} \lambda_k
\Theta_{k,\frac{2m+1}{2}}(\tau)+$$
$$+\frac{\goth{f}(\tau)}{\eta(\tau)}(- \tau) \sum_{j=1}^{2m}
\nu_{j} (\partial \Theta)_{j,\frac{2m+1}{2}}(\tau),$$ for some
constants $\lambda_k$ and $\nu_{j}$. Because of
$$\Theta_{j,k}=\Theta_{-j,k}=\Theta_{2k-j,k}=\Theta_{2k+j,k},$$
$$( \partial \Theta )_{j,k}=-( \partial \Theta )_{-j,k}$$
the previous formula indicates that \be \label{gen-ns-char} \tau
\frac{\goth{f}(\tau)}{\eta(\tau)} (\partial
\Theta)_{j,\frac{2m+1}{2}}(\tau), \ \  j=1,...,m \ee have to be
added to the vector space spanned by irreducible $\striplet$
characters in order to preserve modular invariance. In the case of
the triplet vertex algebra expressions similar to
(\ref{gen-ns-char}) could be interpreted as Miyamoto's
pseudocharacters (cf. \cite{AdM-triplet}). On the other hand, the
$T$ transformation $ \tau \mapsto \tau+1$, maps characters to
supercharacters (multiplied with appropriate scalars). In order to
find an $SL(2,\mathbb{Z})$-closure, we would have to apply the $S$
transformation on the space of supercharacters, but this requires
a knowledge of irreducible $\sigma$-twisted characters. Since we
do not study $\sigma$-twisted $\striplet$-modules in this paper,
at this point we record modular invariance property for the
untwisted sector only.
\begin{theorem} \label{modular-inv} The vector space $\mathcal{NS}$ spanned by: \bea && \chi_{S \Lambda(m+1)}(\tau), \ \chi_{S
\Lambda(i+1)}(\tau), \ \chi_{S
\Pi(m-i)}(\tau), \ \ i=0,...,m-1, \nonumber \\
&& \tau \frac{\goth{f}(\tau)}{\eta(\tau)}( \partial
\Theta)_{m-i,\frac{2m+1}{2}}(\tau), \ \ i=0,...,m-1 \eea is
$(3m+1)$-dimensional and invariant under the subgroup
$\Gamma_{\theta} \subset SL(2,\mathbb{Z})$, where
$\Gamma_\theta=\langle S,T^2 \rangle$.

\end{theorem}

\begin{remark} \label{ramond}
{\em We expect that $S$-transforms of (generalized)
supercharacters are expressible in terms of characters and
generalized characters of $\sigma$-twisted
$\mathcal{SW}(m)$-modules. More precisely, appropriately defined
vector space spanned by characters and generalized supercharacters
, denoted by $\widetilde{\mathcal{NS}}$, and the vector space
spanned by characters and generalized characters of
$\sigma$-twisted modules, denoted by $\mathcal{R}$, should be
inter-related as on the diagram
$$
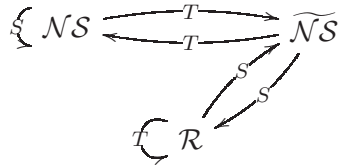
\xymatrix{ \mathcal{NS} \ar@(ul,dl)[]|{S}  \ar@/^/[rr]|{T} &  & \widetilde{\mathcal{NS}} \ar@/^/[ll]|{T} \ar@/^/[dl]|{S} \\ &
\mathcal{R} \ar@(ul,dl)[]|{T} \ar@/^/[ur]|{S} & }$$ It is known
that (super)characters of $N=1$ minimal models in NS and R sector
transform according to this picture (see \cite{IK2}). }
\end{remark}

\section{$\mathcal{SW}(m)$-characters and $q$-series identities}

In this section we discuss fermionic  expressions for irreducible
characters of $\mathcal{SW}(m)$-modules. As we shall see irreducible
$\mathcal{SW}(m)$-modules admit $q$-series formulas similar to those
for the triplet, conjectured by Flohr-Grabov-Koehn \cite{FGK}, and
proven by Warnaar \cite{Wa} ( Feigin et al. independently obtained
similar identities by using different methods \cite{FFT}). More
precisely, the characters of irreducible modules for the super
triplet $\mathcal{SW}(m)$ are intimately related to characters of
irreducible $\mathcal{W}(2m+1)$-modules. It is not clear whether a
deeper connection persists beyond characters.

\subsection{The $m=1$ case: first computation}

Motivated by computations in \cite{FGK} for $\mathcal{W}(2)$, here
we probe double-sum fermionic expressions of irreducible
characters of $\mathcal{SW}(1)$-modules.

As usual, we will be using
$$(a;q)_n=(1-a)(1-aq) \cdots (1-aq^{n-1}),$$
$$(a;q)_\infty=\prod_{i=1}^\infty (1-aq^{i-1}),$$
and sometimes we shall write
$$(q)_n=(q;q)_n,$$
for simplicity.

We start a basic relation \be \label{basic-q}
\displaystyle{ \frac{\prod_{n=1}^\infty
(1+q^{n-1/2})}{\prod_{n=1}^\infty(1-q^n)}=\frac{1}{\prod_{n=1}^\infty
(1-q^{n/2})(1+q^n)}}. \ee
We shall also use a Durfee rectangle identities which
hold for every $k \in \mathbb{Z}_{\geq 0}$, \bea \label{durfee}
&& \frac{1}{\prod_{n \geq 1} (1-q^{n/2})}=\sum_{n=0}^\infty
\frac{q^{(n^2+kn)/2}}{(q^{1/2})_n
(q^{1/2})_{n+k}} \nonumber \\
&& =\sum_{n=0}^\infty \frac{(-q^{1/2})_n
(-q^{1/2})_{n+k}q^{(n^2+kn)/2}}{(q)_n (q)_{n+k}}. \eea Another
useful elementary formula due to Euler is \be \label{euler}
\eta(q)=q^{1/24}\sum_{n=0}^\infty \frac{(-1)^n
q^{(n+1)n/2}}{(q)_n}. \ee

For $m=1$ there are three irreducible characters. We will focus
here on
\be \label{lin-comb-1} \chi_{S
\Lambda(1)}(\tau)=\frac{\goth{f}(\tau)}{\eta(\tau)}\left(\frac{1}{3}
\Theta_{1,\frac{3}{2}}(\tau)+\frac{2}{3}(\partial
\Theta)_{1,\frac{3}{2}}(\tau) \right). \ee

We first notice a theta-function identity
$$(\partial \Theta)_{1,3/2}(\tau)=\frac{\eta(\tau)^3}{\goth{f}(\tau)^2},$$
(essentially, a consequence of the Jacobi triple product identity)
or equivalently
$$\frac{\goth{f}(\tau)}{\eta(\tau)}
(\partial
\Theta)_{1,3/2}(\tau)=\frac{\eta(\tau)^2}{\goth{f}(\tau)}.$$ Now,
we apply the relation
$$\goth{f}(\tau)=\frac{\eta(\tau)^2}{\eta(\tau/2)\eta(2 \tau)}$$
and (\ref{euler}), so we obtain
\bea \label{modify-char} && \frac{\goth{f}(\tau)}{\eta(\tau)}
(\partial \Theta)_{1,3/2}(\tau)=
\eta(2 \tau)\eta(\tau/2)=q^{5/48}
\prod_{n=1}^\infty (1-q^{2n})(1-q^{n/2}) \\
&& = q^{5/48} \sum_{(m_1,m_2) \in \mathbb{Z}_{\geq 0}^2}
\frac{(-1)^{m_1+m_2}
(-q^{1/2};q^{1/2})_{m_2}q^{m_1(m_1+1)+m_2(m_2+1)/4}}{(q^2)_{m_1}(q)_{m_2}}
\nonumber. \eea

On the other hand Durfee square identity (\ref{durfee}) yields
(after some computation)
\bea \label{form-1-2-2} &&
\frac{\goth{f}(\tau)}{\eta(\tau)}\Theta_{1,3/2}(\tau) \nonumber
\\
&& =\frac{q^{5/48}}{(-q;q)_\infty}
\sum_{{\tiny \begin{array}{c} (m_1,m_2) \in \mathbb{Z}_{\geq 0}^2 \\
m_1 \equiv m_2 \ (2) \end{array}}} \frac{(-q^{1/2} ;
q^{1/2})_{m_1}(-q^{1/2} ; q^{1/2})_{m_2}
q^{\frac{3(m_1-m_2)^2}{8}+\frac{(m_1-m_2)}{2}+\frac{m_1
m_2}{2}}}{(q)_{m_1}(q)_{m_2}}.  \eea

Evidently, double fermionic expressions for $(\partial
\Theta)_{1,3/2}(\tau)$ and $\Theta_{1,3/2}(\tau)$ (cf. formulas
(\ref{modify-char}) and (\ref{form-1-2-2}), respectively) appear to have
little in common, so it is unclear to us that
(\ref{lin-comb-1}) admits representation as a closed {\em double}
fermionic sum. Thus, it appears that the $m=1$ case is rather
different compared to the triplet $\mathcal{W}(2)$. This is
perhaps reflected by the fact that $p=2$ triplet admits a
fermionic construction, while such a realization seems to be
absent for $\mathcal{SW}(1)$ and its modules.

\subsection{Irreducible $\mathcal{SW}(m)$ characters from $\mathcal{W}(2m+1)$ characters}

In this part we will be using character formulas of  irreducible $\mathcal{W}(p)$-modules
(see for instance (6.34) and (6.35) in  \cite{AdM-triplet}, or \cite{FHST}).
Recall
$$\goth{f}_2(\tau)=q^{1/24} \prod_{n=1}^{\infty}(1+q^n).$$

The first result in this part is

\begin{proposition} \label{SWW}
\begin{itemize}
\item[(i)] For  $0 \leq i \leq m$, we have
$$\chi_{S\Lambda(i+1)}(\tau)=\frac{\chi_{\Lambda(2i+1)}(\frac{\tau}{2})}{\goth{f}_2(\tau)}.$$
\item[(ii)] For $0 \leq i \leq m-1$, we also have
$$\chi_{S\Pi(m-i)}(\tau)=\frac{\chi_{\Pi(2m-2i)}(\frac{\tau}{2})}{\goth{f}_2(\tau)}.$$
\end{itemize}
Here $\Lambda(i)$ and $\Pi(2m+2-i)$, $i=1, \dots, 2m+1$, are
irreducible $\mathcal{W}(2m+1)$-modules \cite{AdM-triplet}.
\end{proposition}
{\em Proof.} The proof follows from character formulas for
irreducible $\mathcal{W}(p)$-modules, Theorem \ref{char-irr-mod},
and the following  transformation formulas
$$\Theta_{2 j, 2m+1}(\tau/2) = \Theta_{j,\tfrac{2m+1}{2}}(\tau),$$
$$(\partial\Theta)_{2 j, 2m+1}(\tau/2) = 2 (\partial\Theta)_{j,\tfrac{2m+1}{2}}(\tau),$$
$$\frac{1}{\eta(\tau/2)} = \frac{\goth{f}(\tau)}{\eta(\tau)}
\goth{f}_2(\tau).  $$ \qed

We recall two multi-sum identities obtain recently by Warnaar
\cite{Wa} (these identities are essentially conjectures from
\cite{FGK}):

\begin{theorem} For $\lambda=0,\dots ,p$ and $\sigma \in \{0,1\}$ we
have
\bea \label{ole-1} && \sum_{{\tiny \begin{array}{c} n_1,...,n_p=0 \\
n_{p-1}+n_p \equiv 0 \ (2) \end{array}}} \frac{q^{\sum_{i,j=1}^p
B_{i,j}n_i n_j+\lambda/2(n_{p-1}-n_p+\sigma)-\sigma
p/4}}{(q;q)_{n_1} \cdots (q;q)_{n_k}} \\ && =
\frac{1}{(q;q)_{\infty}} \sum_{n \in \mathbb{Z}}
q^{pn^2+(\lambda-\sigma p)n} \nonumber \eea
and
\bea \label{ole-2} && \sum_{{\tiny \begin{array}{c} n_1,...,n_p=0,
\\  n_{p-1}+n_p \equiv 0 \ (2) \end{array}}} \frac{q^{\sum_{i,j=1}^p
B_{i,j}n_in_j+\lambda/2(n_{p-1}+n_p+\sigma)+\sum_{p-\lambda}^{p-2}(i-p+\lambda+1)n_i-\sigma
p/4}}{(q;q)_{n_1} \cdots (q;q)_{n_k}} \\ &&
=\frac{1}{(q;q)_{\infty}} \sum_{n \in \mathbb{Z}}(2n-\sigma+1)
q^{pn^2+(\lambda-\sigma p)n}, \nonumber \eea

where $B_{i,j}$ are entries of the inverse Cartan matrix of the
Lie algebra $D_p$.

\end{theorem}

Equipped with Warnaar's formulas and Proposition \ref{SWW} it is now not hard to prove the
next result
\begin{theorem}

We have the following formulas for irreducible
$\striplet$-characters:

\bea \label{adm-1} && q^{-1/16} \chi_{S
\Lambda(m+1)}(\tau) \\
&&= \sum_{{\tiny \begin{array}{c} n_1,...,n_{2m+1}=0, \\
n_{2m}+n_{2m+1} \equiv 0 \ (2) \end{array}}}
\frac{(-q^{1/2};q^{1/2})_{n_1} \cdots
(-q^{1/2};q^{1/2})_{n_{2m+1}} q^{\sum_{k,l=1}^{2m+1} B_{k,l}n_k
n_l/2}}{(-q;q)_\infty (q;q)_{n_1} \cdots (q;q)_{n_{2m+1}}}. \nonumber \eea

For $i=0,...,m-1$, we have {\tiny \bea \label{adm-3} &&
q^{-a_{i,m}} \chi_{S \Lambda(i+1)}(\tau) \nonumber \\
&& =\sum_{{\tiny \begin{array}{c} n_1,...,n_{2m+1}=0 \\
n_{2m}+n_{2m+1} \equiv 0 \ (2) \end{array}}} \frac{{\tiny
(-q^{1/2};q^{1/2})_{n_1} \cdots (-q^{1/2} ; q^{1/2})_{n_{2m+1}} }
q^{\sum_{k,l=1}^{2m+1}
B_{k,l}n_kn_l/2+(m-i)(n_{2m}+n_{2m+1})/2+\sum_{k=2i+1}^{2m-1}(k-2i)
n_k/2}}{(-q;q)_\infty (q;q)_{n_1} \cdots (q;q)_{n_{2m+1}}},  \nonumber \eea }
and {\tiny \bea \label{adm-4} &&
q^{-b_{i,m}} \chi_{S \Pi(m-i)}(\tau) \nonumber \\
&& =\sum_{{\tiny \begin{array}{c} n_1,...,n_{2m+1}=0 \\
n_{2m}+n_{2m+1} \equiv 1 \ (2) \end{array}}} \frac{{\tiny
(-q^{1/2};q^{1/2})_{n_1} \cdots (-q^{1/2};q^{1/2})_{n_{2m+1}} }
q^{\sum_{k,l=1}^{2m+1} B_{k,l}n_k n_l/2+(m-i)(n_{2m}+n_{2m+1})/2 +
\sum_{k=2i+1}^{2m-1}(k-2i)n_k/2}}{(-q;q)_\infty (q;q)_{n_1} \cdots
(q;q)_{n_{2m+1}}}, \nonumber \eea } where $a_{i,m}$ and $b_{i,m}$
are certain rational numbers.
\end{theorem}
\noindent {\em Proof.} We prove the middle formula only. The other
two formulas follow along the same lines.

Recall that \be \label{super-ferm-1} \chi_{S
\Lambda(i+1)}(\tau)=\frac{\goth{f}(\tau)}{\eta(\tau)}\left(\frac{2i+1}{2m+1}
\Theta_{m-i,\frac{2m+1}{2}}(\tau)+\frac{2}{2m+1}(\partial
\Theta)_{m-i,\frac{2m+1}{2}}(\tau) \right). \ee

Now,
$$\frac{2i+1}{2m+1}
\Theta_{m-i,\frac{2m+1}{2}}(\tau)+\frac{2}{2m+1}(\partial
\Theta)_{m-i,\frac{2m+1}{2}}(\tau)$$
$$={q^{(m-i)^2/(2(2m+1))} \sum_{n \in \mathbb{Z}}
(2n+1)q^{\frac{(2m+1)n^2+2(m-i)n}{2}}}.$$

Finally, if we substitute $q^{1/2}$ for $q$ in (\ref{ole-2}), and
let $p=2m+1$, $\sigma=0$, $\lambda=2m-2i$, and apply formula
(\ref{basic-q}) and simple identity
$$\frac{1}{(q^{1/2};q^{1/2})_{n}}=\frac{(-q^{1/2};q^{1/2})_n}{(q)_n}.$$
the proof automatically follows. \qed

\section{A conjectural relation of $\striplet$ with quantum groups}

Let $\hat{\goth{g}}$ be an untwisted affine Kac-Moody Lie algebra.
Then there is a well-known (Kazhdan-Lusztig) equivalence between
the tensor category of $L_{\goth{g}}(k,0)$-modules $k \in \mathbb{N}$, and the semisimple part of the tensor category
of $U_q(\goth{g})$-modules where $q$ is a certain root of unity (not to be confused with $q=e^{2 \pi i \tau}$ used in the previous section)
depending on the level $k$ and $\goth{g}$ \cite{Fi}. Notice that
on the quantum group side we have a semisimplified category, and
not the full category of $U_q(\goth{g})$-modules.

In \cite{FGST} and \cite{FGST1} (see also \cite{Se}) the authors
proposed a remarkable equivalence between the (enhanced) tensor
category of $\mathcal{W}(p)$-modules and the category of
$\overline{U_q(sl_2)}$-modules, $q=e^{i \pi/p}$, where
$\overline{U_q(sl_2)}$ is the {\em restricted} finite-dimensional
quantum group. While this is still a conjecture for $p>2$, the same
authors established an important weaker equivalence among the
$SL(2,\mathbb{Z})$-module $\goth{Z}_{cft}$ formed by generalized
$\mathcal{W}(p)$ characters and the $SL(2,\mathbb{Z})$-module
$\goth{Z}$, the center of $\overline{U_q(sl_2)}$. Thus, it is
natural question to find Kazhdan-Lusztig dual of the category or
ordinary and logarithmic $\striplet$-modules. In our case the
relevant space of generalized characters is the $\Gamma_\theta$
invariant subspace described in Theorem \ref{modular-inv}, which is
$3m+1$-dimensional.

As indicated in the introduction, we believe that
the quantum group  $U^{small}_q(sl_2)$, $q=e^{\frac{2 i \pi}{2m+1}}$ is relevant for the supertriplet $\striplet$. Here are some evidences.
Firstly, both $\striplet$ and $U^{small}_q(sl_2)$ have the same
number of inequivalent irreducible representations. Also, in
\cite{Ker} (see also \cite{La}) it was proven that the center of
$U^{small}_q(sl_2)$ is $3m+1$-dimensional, and that it carries a
projective action of the modular group. Notice that $3m+1$ is also
(conjecturally) the dimension of the center of $A(\striplet)$.
Thus, in parallel with \cite{FGST}, we expect the following
conjecture to be true.
\begin{conjecture} \label{last-strip} The category of weak $\striplet$-modules is
equivalent to the category of ${{U}^{small}_q(sl_2)}$-modules, where
$q=e^{\frac{2 \pi i}{2m+1}}$.
\end{conjecture}

Finally, Proposition \ref{SWW} is a strong indication for a possible that the
category of $\mathcal{SW}(m)$-modules should be related to a
subcategory of $\mathcal{W}(2m+1)$ and
$\overline{U_q(sl_2)}$-modules, $q=e^{\pi i/(2m+1)}$.

\section{Outlook and final remarks}

There are several research directions we plan to pursue in the future. Let us mention only a few we found
the most interesting.

\begin{itemize}

\item[(i)] The most important problem that we left open
is the existence and description of logarithmic
$\striplet$-modules. We strongly believe the ideas based on modular
invariance as in \cite{AdM-triplet} could be successfully applied for the
super triplet.

\item[(ii)] As with any $N=1$ vertex operator superalgebra, the
most obvious next step would be to examine the category of
$\sigma$-twisted $\striplet$-modules, where $\sigma$ is the parity
automorphisms. As we already indicated (cf. Theorem
\ref{modular-inv}) the space of $SL(2,\mathbb{Z})$-transforms of
irreducible $\striplet$-modules should close a finite-dimensional
vector space. Supposedly characters of irreducible
$\sigma$-twisted modules are included in the same vector space
(cf. Remark \ref{ramond}).

\item[(iii)] Singular vectors in Feigin-Fuchs modules for
the $N=1$ Neveu-Schwarz algebra certainly deserve more attention. We expect
these vectors to have description in terms of modified Jack
polynomials and as kernels of super Calogero-Sutherland
operators. Similar results for the Virasoro algebra have been
obtained in \cite{MY}.

\item[(iv)] Our fermionic expressions for the $\striplet$-characters   indicate a possibility of
parafermionic (or quasiparticle) bases for $\striplet$-modules. For the triplet
$\mathcal{W}(p)$ this problem has been resolved in \cite{FFT}.

\end{itemize}

\section{Appendix}

Here we prove Theorem \ref{2m1prime} and give a strong evidence that  in
Proposition \ref{rel-zhu-prop} the polynomial  $g(x)$ is nonzero for every $m$.  In the process of proving these results
we discovered certain constant term identities which are of independent interest.

We recall
$$U^{F,E}:={\rm Res}_z \frac{(1+z)^{2m}}{z} Y(F,z)E \in \overline{SM(1)}.$$
Then we have
$$U^{F,E}(0):=o(U^{F,E})=\sum_{i \geq 0}{2m \choose i}o( F_{i-1} E) . $$

In Proposition \ref{rel-zhu-druga}  we proved that inside
$A(\overline{SM(1)})$ we have the relation
$$[U^{F,E}]=g([\omega]) [\widehat{H}].$$

Because of the homomorphism from $A(\overline{SM(1)})$ to
$A(\mathcal{SW}(m))$  and Proposition \ref{rel-zhu-prop}  it is
sufficient to show that $U^{F,E}(0)$ acts nontrivially on the top
components of at least one  $SM(1)$-modules $M(1,\lambda) \otimes
F$.

\begin{proposition} \label{gt} Let $v_{\lambda}$ be the highest weight vector in $M(1,\lambda) \otimes F$.
Then we have
$$U^{F,E}(0) \cdot v_{\lambda}=-G_m(t) v_{\lambda},$$
where $t=\langle \lambda,\alpha \rangle$ and
$$
G_m(t)=\sum_{l =1} ^{2m+1} \sum_{i=0}^{l-1} \ \sum_{j=0}^{2m+1 + i
-l} \  \sum_{k = 0}^{l-1-i}
 (-1)^{j+k+l} {2m+1 \choose l}{-2m-1 \choose j} \cdot $$
$$ {-2m-1 \choose k}{2m-t \choose j+k+2m+1}{t \choose i-j-l+2m+1}{t  \choose l-k-1-i}.$$
\end{proposition}
\noindent {\em Proof.}  It is not hard to see that
$$U^{F,E}=\left(\sum_{i=0}^{2m} {\rm Res}_{z_1}{\rm Res}_{z_2}{\rm Res}_{z_3} \frac{(1+z_1)^{2m}}{z_1}z_2^{-i-1}z_3^iY(e^{-\alpha},z_1)Y(e^{\alpha},z_2)
Y(e^{\alpha},z_3)e^{-\alpha} \right)+w,$$ where $w \in
\mathcal{T}$. By repeatedly using the well-known formula (cf.
\cite{LL})
$$E^{+}(\delta,x)E^{-}({\gamma},y)= (1-y/x)^{\langle \delta,\gamma \rangle}E^{-}({\gamma},y)E^+({\delta},x),$$
which holds for every $\delta , \gamma \in \mathbb{Z} \beta$,
we get
$$U^{F,E}=\sum_{i=0}^{2m} {\rm Res}_{z_1}{\rm Res}_{z_2}{\rm Res}_{z_3} \frac{(1+z_1)^{2m}}{z_1}z_2^{-i-1}z_3^i (z_1 z_2 z_3)^{-2m-1} \cdot $$
$$(1-z_2/z_1)^{-2m-1}(1-z_3/z_1)^{-2m-1}(z_2-z_3)^{2m+1}E^{-}(\alpha,z_1)E^{-}(-\alpha,z_2)E^{-}(-\alpha,z_3)+w.$$
Previous formula together with
$$o(E^{-}(\beta,x)) \cdot v_{\lambda}=(1+x)^{-\langle \beta,\lambda \rangle} v_{\lambda}$$
and $$o(w) v_{\lambda}=0$$ implies
$$U^{F,E}(0) \cdot v_{\lambda}=\sum_{i=0}^{2m} {\rm Res}_{z_1}{\rm Res}_{z_2}{\rm Res}_{z_3} \frac{(1+z_1)^{2m}}{z_1}z_2^{-i-1}z_3^i (z_1z_2z_3)^{-2m-1} \cdot $$
$$(1-z_2/z_1)^{-2m-1}(1-z_3/z_1)^{-2m-1}(z_2-z_3)^{2m+1}(1+ z_1)^{-t}(1+z_2)^t (1+z_3)^t v_{\lambda}.$$
The rest follows by expanding generalized rational functions with respect to standard conventions in vertex algebra theory and extracting the residues in all
three variables.
\qed

If we view parameter $t$ as a variable, the expression $G_m(t)$ is a polynomial in $t$ of degree at most $4m+1$.
However, it is a priori not clear that the polynomial $G_m(t)$ is nonzero. We made some computations for small $m$ and we came up with the following hypothesis.

\begin{conjecture} \label{ctc}
$$G_m(t)={2m \choose m}^2 {t+m \choose 4m+1}.$$
\end{conjecture}
We checked this conjecture by using Mathematica package for every $m \leq 20$.

As in Section 11, by using representation theory of $\mathcal{SW}(m)$ it is not hard to see that ${t+m \choose 4m+1}$ must  divide $G_m(t)$ for every $m$.
Since ${\rm deg}(G_m(t)) \leq 4m+1$, then we
have \be \label{am} G_m(t)=A_m {t +m \choose 4m+1 },  \ee  for some constant $A_m$. But even proving $A_m \neq 0$ seems to be a nontrivial problem.

\begin{proposition} \label{3m1} Let $2m+1$ be prime. Then $G_m(t) \neq 0$.
\end{proposition}
{\em Proof.} We will prove this result by virtue of reduction mod $2m+1$. Let $$p=2m+1$$ be a prime. It is not hard to see that in fact  $G_m(a) \in \mathbb{Z}_{(p)}$, for every $a \in \mathbb{Z}$  (in other words, $G_m(a)$ is $p$-integral). Thus it is sufficient to prove that for some $t=t_0$ we have $G_m(t_0) \neq 0 \ {\rm mod} \  p$.
We take $t_0=3m+1$ and examine
$$
G_m(t)=\sum_{l =1} ^{2m+1} \sum_{i=0}^{l-1} \ \sum_{j=0}^{2m+1 + i
-l} \  \sum_{k = 0}^{l-1-i}
 (-1)^{j+k+l} {2m+1 \choose l}{-2m-1 \choose j} \cdot $$
$$ {-2m-1 \choose k}{-m-1 \choose j+k+2m+1}{3m+1 \choose i-j-l+2m+1}{3m+1 \choose l-k-1-i}.$$
The finite sum $G_m(3m+1)$ has many terms divisible by $p$.
 For instance, in the summation, all terms ${ 2m+1 \choose l } \equiv 0 \  {\rm mod}  \ p$ unless
 $l=2m+1$.
After some analysis it is not hard to see that
 possible nontrivial $({\rm mod} \ p)$ contribution comes only if
$k=j=0$ and $l=2m+1$ (in other cases at least one binomial
coefficient is divisible by $p$). Thus we get:

$$G_m(3m+1) \equiv \sum_{i=0}^{2m} (-1) {-m-1 \choose 2m+1}{3m+1 \choose i}{3m+1 \choose 2m-i}  \  {\rm mod}  \  p.$$
Observe the basic relation
$${-m-1 \choose 2m+1}=-{3m+1 \choose 2m+1}=-{3m+1 \choose m}.$$
Also, for $i$ as in the summation we have
$${3m+1 \choose i}{3m+1 \choose 2m-i}  \equiv 0  \  {\rm mod }  \  p,  \  \   i \neq m. $$
However for $i=m$ we have
$${3m+1 \choose m} \equiv 11^{-1} 2 2^{-1} \cdots m m^{-1} \equiv 1 \  {\rm mod}  \  p.$$
Consequently, the summation reduces to a single term
$$G_m(3m+1) \equiv {3m+1 \choose m}^3 \equiv 1 \  {\rm mod} \ p.$$
\qed

Notice that the previous computations support our Conjecture \ref{ctc} because
$${2m \choose m } \equiv  \pm 1  \ {\rm mod}  \  p,$$
so that for $t=3m+1$
$${2m \choose m}^2 {t+m \choose 4m+1}={2m \choose m}^2 {4m+1  \choose 4m+1} \equiv 1 \ {\rm mod }  \  p.$$

\begin{remark} {\em  Because of interesting
arithmetics involved in Propositions \ref{gt} and \ref{3m1}, we
plan to return to Conjecture \ref{ctc} in our future work.}
\end{remark}

\end{document}